\documentclass[12pt, reqno]{amsart}
\usepackage{amssymb}
\usepackage{graphicx}
\usepackage{xcolor} 
\usepackage{tensor}
\usepackage[color=yellow]{todonotes}

\usepackage{enumitem}
\usepackage{fullpage} 
\usepackage{hyperref}

\definecolor{green}{rgb}{0,0.8,0} 

\newtheorem{theorem}{Theorem}[section]

\newtheorem{proposition}[theorem]{Proposition}
\theoremstyle{definition}
\newtheorem{definition}[theorem]{Definition}

\theoremstyle{remark}
\newtheorem{remark}[theorem]{Remark}
\numberwithin{equation}{section}

\newcommand{\tA}{{\tilde A}}

\newcommand{\ta}{{\tilde a}}

\newcommand{\tF}{{\tilde F}}

\renewcommand{\H}{{\mathcal H}}

\newcommand{\nrm}[1]{\Vert#1\Vert}
\newcommand{\abs}[1]{\vert#1\vert}
\newcommand{\brk}[1]{\langle#1\rangle}
\newcommand{\set}[1]{\{#1\}}

\newcommand{\tr}{\textrm{tr}}

\newcommand{\aleq}{\lesssim}
\newcommand{\ageq}{\gtrsim}

\newcommand{\lap}{\Delta}

\newcommand{\ud}{d}
\newcommand{\rd}{\partial}
\newcommand{\nb}{\nabla}

\newcommand{\bb}{\Big}

\newcommand{\alp}{\alpha}
\newcommand{\bt}{\beta}
\newcommand{\gmm}{\gamma}

\newcommand{\dlt}{\delta}

\newcommand{\eps}{\epsilon}
\newcommand{\veps}{\varepsilon}
\newcommand{\kpp}{\kappa}
\newcommand{\lmb}{\lambda}

\newcommand{\sgm}{\sigma}

\newcommand{\bfw}{{\bf w}}

\newcommand{\bfA}{{\bf A}}

\newcommand{\bfD}{{\bf D}}

\newcommand{\bfH}{{\bf H}}

\newcommand{\bfO}{{\bf O}}

\newcommand{\bfQ}{{\bf Q}}

\newcommand{\bbD}{\mathbb D}

\newcommand{\bbH}{\mathbb H}

\newcommand{\bbR}{\mathbb R}
\newcommand{\R}{\mathbb R}
\newcommand{\bbS}{\mathbb S}

\newcommand{\bbZ}{\mathbb Z}

\newcommand{\calC}{\mathcal C}
\newcommand{\calD}{\mathcal D}
\newcommand{\calE}{\mathcal E}
\newcommand{\calF}{\mathcal F}
\newcommand{\calG}{\mathcal G}
\newcommand{\calH}{\mathcal H}

\newcommand{\calO}{\mathcal O}
\newcommand{\calP}{\mathcal P}
\newcommand{\calQ}{\mathcal Q}

\setlist[enumerate]{leftmargin=2em, label=(\arabic*)}
\setlist[itemize]{leftmargin=2em}

\newcommand{\covD}{\bfD}

\newcommand{\Diff}{\mathrm{Diff}}

\newcommand{\Cal}{\text{Cal}}

\newcommand\DA{{\mathbf{DA}}}
\newcommand\DB{{\mathbf{DB}}}

\newcommand{\g}{\mathfrak  g}
\newcommand{\G}{\mathbf{G}}
\newcommand{\la}{\langle}
\newcommand{\ra}{\rangle}
\newcommand{\p}{\partial}
\renewcommand{\aa}{\alpha}
\renewcommand{\bb}{\beta}
\renewcommand{\P}{\mathbf P}
\newcommand{\rst}{\!\upharpoonright}		

\newcommand{\anb}{\displaystyle{\hskip-.1em \not \hskip-.25em \nb}}

\newcommand{\nE}{{\calE}} 
\newcommand{\spE}{\calE_{e}} 
\newcommand{\Egs}{E_{\text{GS}}} 

\newcommand{\hM}{\calQ} 
\newcommand{\ch}{\boldsymbol{\chi}}

\begin{document}

\title[]{The Threshold Theorem for  the $(4+1)$-dimensional
Yang--Mills  equation: an overview of the proof }
\author{Sung-Jin Oh}%
\address{KIAS, Seoul, Korea 02455}%
\email{sjoh@kias.re.kr}%

\author{Daniel Tataru}%
\address{Department of Mathematics, UC Berkeley, Berkeley, CA, 94720}%
\email{tataru@math.berkeley.edu}%


\date{\today}%
\begin{abstract}
  This article is devoted to the energy critical hyperbolic
  Yang--Mills system in the $(4+1)$ dimensional Minkowski space, which
  is considered by the authors in a sequence of four papers \cite{OTYM1},
  \cite{OTYM2}, \cite{OTYM2.5} and \cite{OTYM3}. The final outcome of
  these papers is twofold: (i) the Threshold Theorem, which asserts
  that global well-posedness and scattering hold for all topologically trivial
  initial data with energy below twice the ground state energy, and 
  (ii) the Dichotomy Theorem, which for larger data in arbitrary topological
  classes provides a choice of two outcomes,
  either a global, scattering solution or a soliton bubbling off.  In
  the last case, the bubbling off phenomena can happen either (a) in
  finite time, triggering a finite time blow-up, or (b) in infinite time.
  Our goal here is to describe these results, and to provide an
  overview of the flow of ideas within their proofs in \cite{OTYM1},
  \cite{OTYM2}, \cite{OTYM2.5} and \cite{OTYM3}. 
\end{abstract}
\maketitle
\setcounter{tocdepth}{1}
\tableofcontents

\section{Introduction}

\subsection{Lie groups and Lie algebras}
Let $\G$ be a compact noncommutative Lie group and $\g$ its associated Lie
algebra. We denote by $Ad(O) X = O X O^{-1}$ the action of $\G$ on
$\g$ by conjugation (i.e., the adjoint action), and by $ad(X) Y = [X,
Y]$ the associated action of $\g$, which is given by the Lie
bracket. We introduce the notation $\brk{X, Y}$ for a bi-invariant inner product on $\g$,
\begin{equation*}
\la [X,Y],Z \ra = \la X, [Y,Z] \ra, \qquad X,Y,Z \in \g, 
\end{equation*}
or equivalently 
\begin{equation*}
\la X,Y \ra = \la Ad(O) X,  Ad(O) Y  \ra, \qquad X,Y \in \g, \quad O \in \G. 
\end{equation*}
If $\G$ is semisimple then one can take 
$\brk{X, Y} = -\tr(ad(X) ad(Y))$ i.e. negative of the Killing form on $\g$, which is then positive definite,
However, a bi-invariant inner product on $\g$ exists for any compact Lie group $G$.

\subsection{The Yang--Mills evolution}
Let $\R^{4+1}$ be the five dimensional Minkowski space  with
the standard Lorentzian metric $m= \text{diag}(-1,1,1,1,1)$.
Denote by $A_{\alpha}: \R^{4+1}\rightarrow \g$, $\alpha = 0, \ldots,4$,
a connection $1$-form taking values in the Lie algebra $\g$, and by $\covD_\alpha$
the associated covariant differentiation,
\[
\covD_{\alpha} B:= \partial_{\alpha} B + [A_{\alpha},B],
\]
acting on $\g$-valued functions $B$.
Introducing the curvature $2$-form
\[
F_{\alpha\beta}: = \partial_{\alpha}A_{\beta}
- \partial_{\beta}A_{\alpha} +[A_\alpha,A_\beta],
\]
the \emph{hyperbolic Yang--Mills equation} is the Euler--Lagrange
equation associated with the formal Lagrangian action functional
\[
\mathcal{L}(A_{\alpha}): = \frac{1}{2}\int_{\R^{4+1}} \la F_{\alpha\beta}, F^{\alpha\beta}\ra \,dxdt.
\]
Here we are using the standard convention of raising indices using the metric $m$. Thus,
the Yang--Mills equation takes the form
\begin{equation}\label{ym}
\covD^\alpha F_{\alpha \beta} = 0.
\end{equation}
There is a natural energy-momentum tensor associated to the Yang--Mills
equations, namely 
\[
T_{\alpha \beta} =    2 m^{\gamma \delta} \la F_{\alpha \gamma}, F_{\beta \delta}\ra
- \frac12 m_{\alpha \beta} \la F_{\gamma \delta}, F^{\gamma \delta}\ra.
\]
If $A$ solves the Yang--Mills equation \eqref{ym} then $T_{\alpha \beta}$ is divergence-free,
\begin{equation}\label{divT}
\partial^\alpha T_{\alpha \beta} = 0.
\end{equation}
Integrating this for $\beta = 0$ yields a conserved energy
\begin{equation}\label{energy}
\nE(A) = \nE_{\set{t} \times \bbR^{4}}(A) = \int_{\set{t} \times \R^4} T_{00}\, dx = \sum_{\alp < \bt} \int_{\set{t} \times \R^4} \brk{F_{\alp \bt}, F_{\alp \bt}} \, d x.
\end{equation}
 The case $\beta \neq 0$ yields further conservation laws, i.e. the momentum,
which play no role in the present work.

 The Yang--Mills equation also has a scale invariance property,
\[
A(t,x) \to \lambda A(\lambda t,\lambda x) \quad (\lmb > 0).
\]
The energy functional $\nE$ is invariant with respect to scaling precisely in dimension
$4+1$. For this reason we call the $4+1$ problem \emph{energy critical};  this is one of the
motivations for our interest in this problem.

\subsection{Gauge invariance} 
In order to study the Yang--Mills equation as a well-defined evolution
in time one needs to also consider its gauge invariance. Given a map $O = O(t,x)$ taking values in the group $\G$,
we introduce
\begin{equation*}
	O_{;\alp} = \rd_{\alp} O O^{-1},
\end{equation*}
which now takes values in the Lie algebra $\g$. The gauge transformation of a connection $A$ by $O$ is
   \[
A_\alpha \longrightarrow Ad(O) A_{\alp} - O_{; \alp} =: \calG(O) A_{\alp},
\]
which makes the associated differentiation $\covD$ covariant with $Ad(O)$.
Correspondingly, the curvature tensor changes by
  \[
F_{\alpha \beta} \longrightarrow Ad(O) F_{\alpha\beta}.
\]
Clearly, the Yang--Mills equation \eqref{ym} is invariant under such transforms.

As a consequence, solutions are a-priori defined as equivalence classes.  In order
to uniquely select representatives for the solutions to the
Yang--Mills equation within each equivalence class one needs to add an
additional set of constraint equations; this procedure is known as
{\em gauge fixing}. This issue is fundamental for the fine analysis of
the Yang--Mills equation.  In choosing a gauge, one is naturally led
to pursue conflicting goals:

\begin{enumerate}
\item[(i)] Causality: the system should have finite speed of propagation
\item[(ii)] Structure: the nonlinearity should exhibit null structure type cancellation
\item[(iii)]  Large data: the gauge should be well-defined for large data.
\end{enumerate}

Historically there are (at least) three gauges that have played a role in the study of 
the hyperbolic Yang--Mills evolution:
\smallskip

\emph{1. The Lorenz gauge}, 
\[
\partial^\alpha A_\alpha = 0.
\]
In this gauge the Yang--Mills equation becomes a system of
semilinear wave equations for $A_\alpha$, and in particular it has
finite speed of propagation.  This gauge is very convenient for local
well-posedness for large but regular data.  However, it is not so good
in the low regularity setting as it does not capture well the null
structure, see e.g. \cite{MR3519539}.
\smallskip

\emph{ 2. The temporal gauge}, 
\[
A_0 = 0.
\]
This again insures that the above system is strictly hyperbolic, and
in particular it has finite speed of propagation.  In this gauge the
equations can be understood as a semilinear wave equation for the curl
of $A_x$, coupled with a transport equation for its divergence. This
gauge is also very convenient for local well-posedness for large but
regular data, and it fully describes all regular
solutions to the hyperbolic Yang--Mills equation.
Again there are multiple technical difficulties if one tries to
implement such a gauge in the low regularity setting or globally in 
time. In particular we have  no dispersion for the divergence of $A$.
This gauge will play an auxiliary role in our analysis, and is described 
in greater detail in Section~\ref{sec:topo}.

\smallskip

\emph{3. The Coulomb gauge,}
\[
  \sum_{j=1}^4\partial_j A_j = 0.
\]
Here the causality is lost; however, the Coulomb gauge is an
``elliptic'' gauge which captures well the null structure of the
problem, and thus works well in low regularity settings. Indeed, the
Coulomb gauge was used in \cite{KT} to prove the small data result for
this problem. Unfortunately, it seems that the Coulomb gauge cannot be
implemented globally for large data, even after restricting to those below the ground state energy.
Nevertheless,  for expository purposes we do provide a brief review of the
Coulomb gauge in the beginning of Section~\ref{sec:caloric}. 

For the reasons described above, these three gauges seem inadequate
for the purpose of proving the Threshold Theorem (to be described below). Instead, in our first article 
\cite{OTYM1} we introduce a new gauge, namely 
\medskip

\emph{4. The caloric gauge}. This  is defined via the Yang--Mills heat
flow and is  described  in Section~\ref{sec:caloric}.  It  has the key property that it
is globally defined for all data below the ground state energy. In
addition, to the leading order this agrees with the Coulomb gauge, so
there are many similarities between the  analysis in the caloric and Coulomb gauges.

\subsection{Yang--Mills initial data sets}
In order to consider the hyperbolic Yang--Mills problem as an evolution equation
we need to consider initial data sets. An \emph{initial data set} for
\eqref{ym} consists of a pair of 1-forms $(a_{j}, e_{j})$  on $\bbR^{4}$. 
We say that $(a_{j},e_{j})$ is the initial data for a Yang--Mills solution $A$ if
\begin{equation*}
	(A_{j}, F_{0 j}) \rst_{\set{t=0}} = (a_{j}, e_{j}).
\end{equation*}
The curvature of $a$ is denoted by $f$ in what follows. 

Note that \eqref{ym} imposes the condition that the following
equation be true for any initial data for \eqref{ym}:
\begin{equation} \label{eq:YMconstraint}
	\covD^{j} e_{j} = 0.
\end{equation}
where $\covD^j$ denotes the covariant derivative with respect to the $a_j$
connection.  This equation is the \emph{Gauss} (or the
\emph{constraint}) \emph{equation} for \eqref{ym}. 

\begin{definition}
\begin{enumerate}
\item A regular initial data set  for the Yang--Mills equation  
is a pair of $1$-forms  $(a_j,e_j) \in H_{loc}^N \times H^{N-1}$, $N \geq 2$, with $f \in H^N$,
and which satisfies the constraint equation \eqref{eq:YMconstraint}.
\item A finite energy  initial data set  for the Yang--Mills equation  
is a pair  of $1$-forms $(a_j,e_j) \in H^1_{loc} \times L^2$ with $f \in L^2$ and
which satisfies the constraint equation \eqref{eq:YMconstraint}.
\end{enumerate}
\end{definition}

\subsection{Yang--Mills solutions}
We begin by defining the notions of regular and finite energy solutions:

\begin{definition}
\begin{enumerate}
\item Let $N \geq 2$. A regular solution  for the Yang--Mills equation in an open set $\calO \subset \R^{4+1}$  
is a connection $A \in C([0,T];H^N_{loc})$, whose curvature satisfies $F \in C([0,T];H^{N-1}_{loc})$
and which solves the equation \eqref{ym}.
\item A finite energy solution  for the Yang--Mills equation in the open set $\calO$
is a connection $A \in C([0,T];H^1_{loc})$, whose curvature satisfies $F \in C([0,T];L^2_{loc})$
and which is the limit of regular solutions in this topology.
\end{enumerate}
\end{definition}

We carefully remark that this definition does not require a gauge choice. Hence at this point
solutions are still given by equivalence classes.  Corresponding to the 
above classes of solutions, we have the classes of gauge transformations
which preserve them:

\begin{definition}
\begin{enumerate}
\item Let $N \geq 2$. A regular gauge transformation in an open set $\calO \subset \bbR^{4+1}$ is 
 is a map 
\[
O: \calO \to \G
\]
with the following regularity properties:
\[
O_{;x}, O_{;t} \in C_t (H^{N+1}_{loc}).
\]
\item An admissible gauge transformation in an open set $\calO \subset
\bbR^{4+1}$ is a similar map with the following regularity properties:
\[
O_{;x}, O_{;t} \in C_t (H^1_{loc}).
\]
\end{enumerate}
\end{definition}

Using this notion we can now talk about gauge equivalent connections:

\begin{definition}
 Two  finite energy connections $A^{(1)}$ and $A^{(2)}$  in an open set $\calO \subset \bbR^{4+1}$ are gauge equivalent
if there exists an admissible gauge transformation $O$ so that $A^{(2)} = O A^{(1)} O^{-1} - O_{;x}$.
\end{definition}

\subsection{Topological classes}
The space of finite energy Yang--Mills connections in $\R^4$ is not
connected. Instead, such connections can be classified in terms of
their \emph{topological class}; see Section~\ref{sec:topo} for more details. 

For a compact base manifold, such as $\bbS^{4}$, this term refers to
the isomorphism classes of principal $\G$-bundles which supports the
connection. On the other hand, for $\bbR^{4}$, which is contractible
and thus supports only the trivial fiber bundles, a topological class
must be interpreted rather as a property of a connection.

 In the particular case of four dimensional $SU(2)$ connections  the topological 
class is easily described in terms of the
(second) Chern number
\[
c_2 = \frac{1}{8 \pi^2} \int_{\R^4} \tr ( F \wedge F).
\]
This is always an integer if $A$ has finite energy. For an arbitrary compact noncommutative
Lie group, we have an analogue of $c_{2}$, 
\[
\ch(A) = \int_{\bbR^{4}} - \brk{F \wedge F} = \frac{1}{4} \int_{\bbR^{4}} - \brk{F_{ij}, F_{k \ell}} \, d x^{i} \wedge d x^{j} \wedge d x^{k} \wedge d x^{\ell},
\]
which we denote by $\ch(A)$ and call the \emph{characteristic number}. This quantity is still a topological invariant, but it no longer fully describes the topological class.

The connections which are in the same class as the zero connection are
called \emph{topologically trivial}.  For such connections, $\ch = 0$. An alternative way to describe
topologically trivial connections is given by the following result,
which generalizes Uhlenbeck's lemma ~\cite{MR648356}:

\begin{theorem}[\cite{OTYM2.5}] A finite energy connection $A$ in $\R^4$ 
is topologically trivial iff $A \in \dot H^1$ in a suitable gauge.
\end{theorem}

A further ``Good Global Gauge Theorem'' is provided in \cite{OTYM2.5} for 
finite energy connections which are not topologically trivial.

\subsection{Solitons and the ground state energy}
Steady states for the hyperbolic Yang--Mills equation are called
\emph{harmonic Yang--Mills connections}, and play an important role in
our work.  They solve the equations
\begin{equation}\label{ym-e}
\covD^j F_{kj} = 0 \qquad \hbox{in } \R^4,
\end{equation}
and can be seen as critical points for the Lagrangian
\[
\spE(A) = \frac{1}{2}\int_{\R^{4+1}} \la F_{ij}, F^{ij}\ra \,dxdt.
\]
The key elliptic regularity result is as follows:

\begin{theorem} [Uhlenbeck \cite{MR648356, U2}]
 $\dot H^1$  harmonic Yang--Mills connections are smooth in a suitable gauge. 
\end{theorem}

The question of existence of finite energy harmonic Yang--Mills 
connections is best phrased in terms of the topological classes  described above:
\begin{theorem} \label{thm:gs}
 The following properties hold for harmonic Yang--Mills connections:
\begin{enumerate}
\item Within each topological class there exist energy
  minimizers. These are called instantons, and come in two varieties,
  self-dual $F = \star F$ and anti-self-dual $F = -\star F$, depending on the topological class. 

\item In particular, there exists a unique (up to symmetries) minimal energy 
 nontrivial harmonic Yang--Mills connection $Q$, which is necessarily an instanton, whose energy $\Egs$ satisfies
\[
\nE(Q_{GS}) = |\ch(Q_{GS})|.
\]
\item All  nontrivial harmonic Yang--Mills connection $a$, with energy $\nE(Q) \leq 2 \Egs$ 
are instantons and satisfy
\[
\nE(Q) = |\ch(Q)|.
\]
\end{enumerate}
\end{theorem}

Parts (1) \& (2) are classical. We remark that part (3), which follows from a recent result of \cite{GKS}, is nontrivial due to existence of non-minimizing harmonic Yang--Mills connections \cite{SSU}. We refer to \cite[Sections~1.8 and 6]{OTYM2.5} for further discussion.

As a consequence of the above properties, it easily follows that in 
the class of topologically trivial connections, the threshold for nontrivial
harmonic Yang--Mills connections is $2 \Egs$ rather than $\Egs$.

We also remark that harmonic Yang--Mills connections which are not
energy minimizers no longer have to be self-dual or anti-self-dual.

The harmonic Yang--Mills connections are relevant for the hyperbolic
Yang--Mills flow for multiple reasons. First of all, they provide
examples of solutions that do not scatter.  Further, above the ground
state energy $\Egs$ there are examples of solutions which blow up in
finite time, with a profile which approaches a rescaled instanton, see
\cite{KST2, MR2929728}. Thus, the ground state energy arises as a natural
threshold in the large data well-posedness theory, and one is led to
the \emph{Threshold Conjecture}, which asserts that the Yang--Mills
problem is globally well-posed below the ground state energy.
All such connections must be topologically trivial. However, as discussed above,
for such connections the correct threshold  is $2 \Egs$. 
Based on the above discussion, we will call 
\emph{subthreshold data/solution} any topologically trivial hyperbolic
Yang--Mills data/solution with energy below $2\Egs$.

\subsection{The main results }

The main question we are concerned with is whether the hyperbolic
Yang--Mills equation \eqref{ym}  is globally well-posed in the space of finite energy connections
 in the $4+1$ dimensional setting.  The small
data global well-posedness was recently proved by Krieger together
with the second author in \cite{KT}, so our main interest here is in
large solutions.  The Threshold Conjecture asserts that global
well-posedness in the energy space holds below the ground state energy. 
 \medskip

The first goal of our four  papers \cite{OTYM1, OTYM2, OTYM2.5, OTYM3} is to
establish the validity of (a more precise form of) this conjecture. In the simplest form, our
result can be phrased as follows:

\begin{theorem}[Threshold Theorem for Energy Critical Yang--Mills] \label{thm:threshold}
Global well-posedness \\ and scattering holds for the energy critical hyperbolic 
Yang--Mills evolution in $\R^{4+1}$ for all topologically trivial initial data with
  energy below $2\Egs$.
\end{theorem}

Since scattering solutions are necessarily topologically trivial, we are justified in considering only the topologically trivial data in Theorem~\ref{thm:threshold}.
This restriction, in view of Theorem~\ref{thm:gs}, is the reason why our threshold is $2\Egs$ rather than just $\Egs$.

The statement of this theorem should be understood as follows:
\begin{itemize}
\item For each smooth subthreshold initial data $(a,e)$ there exists a global
  smooth solution, which is unique up to gauge transformations.

\item For each subthreshold data in $\dot{H}^1 \times L^2$   there exists a
  solution $(A, \rd_{t} A) \in C(\R ; \dot H^1 \times L^2)$ which is the unique
  limit of smooth solutions up to gauge transformations.
\end{itemize}

The above formulation of the result is gauge independent. However, in
order to both prove this result and to provide a better description of
the solutions, including their scattering properties, it is essential
to fix the gauge choice in a favorable way.  For our problem, the
classical choices of gauge (Lorenz, temporal or Coulomb) seem to present 
different but equally insurmountable difficulties. We instead rely on the \emph{caloric
  gauge}, which is constructed based on the regularity theory of the
Yang--Mills heat flow, the parabolic counterpart of \eqref{ym}.  A gauge dependent 
formulation of this result will be provided later on, see Theorem~\ref{t:main}.

\medskip

The second goal of our four papers \cite{OTYM1, OTYM2, OTYM2.5, OTYM3}
is to also consider solutions which do not satisfy the topological and
energy constraint of the Threshold Theorem. Then on the one hand, we know there exist
solutions which blow-up or are global but do not scatter, see \cite{KST2,MR2929728}.  On the other hand
scattering can only hold for topologically trivial solutions.  Because
of this, our second result offers a dichotomy:

\begin{theorem} [Dichotomy Theorem for Energy Critical Yang--Mills]\label{thm:dichotomy}\
\!\!\!\!\! The energy critical hyperbolic Yang--Mills evolution in $\R^{4+1}$
  is locally well-posed in the energy space. Further, one of the following two properties
  must hold for the maximal solution:

\begin{enumerate}[label=(\roman*)] 
\item The solution is topologically trivial, global and scatters at infinity. 

\item The solution bubbles off a soliton either
\begin{enumerate}[label=(\alph*)] 
\item at a finite blow-up time, or

\item at infinity.
\end{enumerate}
\end{enumerate}

\end{theorem}
We note that these two alternatives hold separately for positive and negative time.
In other words we do not eliminate the scenario where, say, scattering holds 
for positive time  while finite time blow-up occurs for negative time.

To fully describe this result we need to clarify the meaning of bubbling off.
We do this in the two scenarios, of finite time blow-up solutions and of global solutions.

\medskip 
\emph{a) The finite time blow-up scenario:} Let $t_0 > 0$ be the
blow-up time (maximal existence time) for a finite energy Yang--Mills
connection $A$. By energy conservation, finite speed of propagation
and the small data result there must exist a point $x_0 \in \R^4$ so
that energy concentrates in the backward blow-up cone centered at $(t_0,x_0)$,
namely $C = \{ |x-x_0| < t_0-t \}$, in the sense that
\[
\lim_{t \nearrow  t_0} \nE_{S_{t}} (A) > 0 .
\]
where $S_{t} = C \cap (\set{t} \times \bbR^{4})$.

In this context, we say that \emph{$A$ bubbles off a soliton  at $(t_0,x_0)$} if 
 there exists a sequence of points $ (t_n,x_n) \to (t_0,x_0)$ and scales $r_n$ with the 
following properties:
\begin{enumerate}
\item Time-like concentration,
\[
\limsup_{n\to \infty} \frac{x_n-x_0}{|t_n-t_0|}  = v, \qquad |v| < 1 
\]

\item Below self-similar scale, 
\[
\limsup_{n\to \infty} \frac{r_n}{|t_n-t_0|} = 0
\]

\item Convergence to soliton:
\[
\lim_{n \to \infty} r_{n} \calG(O_n) A(t_{n} + r_{n} t, x_{n} + r_{n} x) = L_v Q(t,x) 
\quad \hbox{in} \ H^1_{loc}([-1/2, 1/2] \times \bbR^{4}) 
\]
for some sequence of admissible gauge transformations $O_n$, a Lorentz transformation
$L_v$ and finite energy harmonic Yang--Mills connection $Q$.
\end{enumerate}
We remark that for a finite energy harmonic Yang--Mills connection $Q$ we must have 
\[
\nE(Q) \leq \nE(L_v Q)
\]
with equality iff $v = 0$.

\medskip

\emph{ b) Global solutions.} Here we consider  a finite energy Yang--Mills
connection $A$ which is global forward in time. We say that $A$ \emph{bubbles off a soliton  at
infinity} if there exists a sequence of points $C \ni (t_n,x_n) \to \infty$ and scales $r_n$ with the 
following properties:
\begin{enumerate}
\item Time-like concentration,
\[
\limsup_{n\to \infty} \frac{x_n}{t_n}  = v, \qquad |v| < 1 
\]

\item Below self-similar scale, 
\[
\limsup_{n\to \infty} \frac{r_n}{t_n} = 0
\]

\item Convergence to soliton:
\[
\lim_{n \to \infty} r_{n} \calG(O_n) A(t_{n} + r_{n} t, x_{n} + r_{n} x) = L_v Q(t,x) 
\quad \hbox{in} \ H^1_{loc}([-1/2, 1/2] \times \bbR^{4}) 
\]
for some sequence of admissible gauge transformations $O_n$, a Lorentz transformation
$L_v$ and finite energy harmonic Yang--Mills connection $Q$.
\end{enumerate}

The proof of these two theorems is the final outcome of the sequence of papers
 \cite{OTYM1},   \cite{OTYM2}, \cite{OTYM2.5} and \cite{OTYM3}.
These contain conceptually disjoint, self-contained logical steps
which address different aspects of the problem, as follows:
\medskip

\begin{description}[itemsep=5pt]
\item  [I. The caloric gauge \cite{OTYM1}] This first paper uses the
  \emph{Yang--Mills heat flow} in order to introduce the \emph{caloric
    gauge}, which is central in our analysis. Its main outcome is to
  provide a complete caloric gauge representation for the hyperbolic
  Yang--Mills equation \eqref{ym}. Along the way, we also establish the
  Threshold and the Dichotomy Theorems for the Yang--Mills heat flow. In particular, the former allows us to prove 
that all subthreshold data admit a caloric representation.  These results are discussed
in Section~\ref{sec:caloric}. 

\item[II. Energy dispersed solutions  \cite{OTYM2}] Here we develop the
analytic tools which are needed in order to understand the hyperbolic Yang--Mills flow
in the caloric gauge. The main result is a strong quantitative a-priori bound 
for \emph{energy dispersed solutions}, which in particular implies local well-posedness
as well as small data global well-posedness in the caloric gauge. The notion
of energy dispersion as well as the main results are described in Section~\ref{sec:ed}.

\item [III. Large data and causality \cite{OTYM2.5}] Since not all
  Yang--Mills solutions can be placed in the caloric gauge, in this
  article we show how to switch the qualitative part of the analysis
  (but not the analytic part) into the temporal gauge, in order to be
  able to deal with data with above threshold energy. The overview in
  Section~\ref{sec:topo} also covers topological classes, initial data surgery
   and gauge matters such as patching of local solutions.  

  \item[ IV. Blow-up analysis \cite{OTYM3}] In this final step we use
    Morawetz type bounds in order to perform a blow-up analysis which
    leads to the proof of the two theorems above.  This is where the
    results in the previous two papers \cite{OTYM2} and \cite{OTYM3}
    are used, but not the the analysis leading to these results. This
    is described in the last section.

\end{description}

We finally remark that these papers build upon a large body of
work. This begins with early results on Yang--Mills above scaling
\cite{MR579231, MR649158, MR649159, KlMa1, Kl-Tat}, where the structure
of the equations was first understood and exploited. Our general
approach broadly follows the outline of similar results for wave maps,
starting with the small data problem, the null frame function spaces
and the renormalization idea \cite{Tat,Tao2,Tat2} and continuing with
the induction on energy based energy dispersion approach in the proof
of the Threshold and Dichotomy Theorem in \cite{ST1, ST2} (see also
\cite{KriSch} and \cite{Tao:2008wo, Tao:2008wn, Tao:2008tz, Tao:2009ta,
  Tao:2009ua}). The similar results for the closely related Maxwell-Klein-Gordon
equation at critical regularity were proved in the small data case in \cite{RT} ($d \geq 6$)
and \cite{KST} ($d \geq 4$), respectively large data in \cite{OT1, OT2, OT3} and independently
in \cite{KL}. Finally, the small data results for (YM) were obtained only recently 
in \cite{KS} ($d \geq 6$) and \cite{KT} ($d \geq 4$). For a more extensive overview of related literature we refer the reader to
\cite{OTYM3}.   Some further comments are provided in each of the following 
sections as needed.

\subsection*{Acknowledgments} 
Part of the work was carried out during the semester program ``New Challenges in PDE''
held at MSRI in Fall 2015. S.-J. Oh was supported by the Miller Research Fellowship from the Miller Institute, UC Berkeley and the TJ Park Science Fellowship from the POSCO TJ Park Foundation. D. Tataru was partially
supported by the NSF grant DMS-1266182 as well as by the Simons
Investigator grant from the Simons Foundation.

\section{The caloric gauge} \label{sec:caloric} This section describes
the main results of \cite{OTYM1}, whose aim is to develop the
\emph{caloric gauge} as our main gauge of choice in the study of the
hyperbolic Yang--Mills evolution.  

Let us take as a starting point of our discussion the following small data result proved earlier  in \cite{KT}:
\begin{theorem} \label{t:coulomb}
The hyperbolic Yang--Mills equation in $\bbR^{4+1}$   is globally well-posed in the Coulomb gauge for 
all initial data with small energy.   
\end{theorem}
 Unfortunately, while  the Coulomb gauge works well in 
the small data problem, it does not appear to work for large data, even after restricting to
only subthreshold data. This large data difficulty with the Coulomb gauge compels us to look
for a different gauge choice, in which the Yang--Mills equation
exhibits a similar null structure as the Coulomb gauge, yet which can
be used in the large data problem. 

 Our solution to this problem is to
introduce and use the \emph{(global) caloric gauge}, which is constructed with
the help of the \emph{Yang--Mills heat flow}.  A more localized form
of this gauge was previously introduced by the first author in
\cite{Oh1, Oh2}, in order to study local well-posedness questions for
the $3+1$ dimensional hyperbolic Yang--Mills equation. This was in turn inspired
by Tao's caloric gauge for wave maps \cite{Tao-caloric}, which is based
on the harmonic map heat flow.

On the one hand, the caloric gauge resembles Coulomb gauge in the
sense that a generalized Coulomb condition holds (to be discussed in
more detail in Section~\ref{subsec:caloric-pf}). On the other hand, it
can be used for a larger class of connections, which in particular
includes all subthreshold connections (essentially by the Threshold
Theorem for the Yang--Mills heat flow, see Theorem~\ref{t:caloric}
below). Therefore, it furnishes a natural setting to state and prove the Threshold Theorem for the hyperbolic Yang--Mills equation; see Theorem~\ref{t:main} below.

\subsection{The Coulomb gauge and the null structure}
\label{subsec:coulomb}
Before we describe the caloric gauge, we first review the null
structure of the hyperbolic Yang--Mills equation in the Coulomb gauge,
which plays essential role in low regularity problems for the
Yang--Mills equation.

Consider the expansion of the Yang--Mills equation \eqref{ym} in terms
of $A$, which takes the form
\begin{equation}\label{eq:ym-exp}
\Box A_\bb + 2 [ A_\aa,\p^\aa A_\bb] = 
 \p_\bb \p^\aa A_\aa -[\p^\aa  A_\aa,A_\bb] 
 + [A^\aa,  \p_\bb A_\aa]  - [A^\aa,[A_\aa,A_\bb]]. 
\end{equation}
where $\Box_{A}: = \covD^{\alpha}\covD_{\alpha}$ is the covariant d'Alembertian (or the covariant wave operator).
Separating the spatial part and the temporal part of the connection,
one immediately sees that the spatial divergence of the solutions 
plays a prominent role. Precisely, one can rewrite the equations in the form
\begin{equation}\label{eq:div1}
\begin{split}
\Box_A A_j  = & \ 
 \p_j  \p^k A_k + \p_j \p^0 A_0   
 + [A^\aa,  \p_j A_\aa]   
\\
\Delta_A A_0  = & \
 \p_0 \p^j A_j  
 + [A^j,  \p_0 A_j]  .
\end{split}
\end{equation}
Thus, when imposing the Coulomb gauge condition, 
\begin{equation}\label{cg}
  \sum_{j=1}^4\partial_j A_j = 0,
\end{equation}
the above equations turn into a hyperbolic system for the main variables
\[
\Box_A A_j  =   \p_j \p^0 A_0    + [A^\aa,  \p_j A_\aa].   
\]
In order to eliminate the first term on the right and also to restrict the evolution 
to divergence free fields $A_j$ we apply the Leray projection $\P$, and rewrite  
the equation in the form
\begin{equation}\label{ym-cg}
\Box A_j  = \P \left(  [A^\aa,  \p_j A_\aa] - 2 [ A^\alpha, \partial_\alpha A_j] 
- [ \partial_0 A_0, A_j]  - [A^\alpha,[A_\alpha,A_j]] \right).  
\end{equation}

Here the $A_0$ component plays an auxiliary role,
and is determined at each fixed time via the elliptic equation
\begin{equation}\label{ym-cg-a0}
\Delta_A A_0  =    [A^j,  \p_0 A_j]  .
\end{equation}
This does not yet yield a self-contained system, as the time derivative of $A_0$ 
also appears in the first equation. A slightly more involved computation yields
the equation
\begin{equation}\label{ym-cg-dta0}
\partial^j \covD_j \covD^0 A_0 =  \partial^j \left( 2 [ A_0,\partial^0 A_j]+   
 [\partial_j A_{\alpha}, A^\alpha] + [A_\alpha,[A^\alpha,A_j] ] \right)      
\end{equation}
which serves to also determine $\covD^0 A_0 $ in an elliptic fashion.

As one can easily see above, the Yang--Mills equations in the Coulomb
gauge can be viewed as an evolution equation \eqref{ym-cg} for the
spatial part $A_x$ of the connection, whereas $A_0$ and $\covD^0 A_0 $
play the role of auxiliary, dependent variables.  All terms in the
equation which involve $A_0$ can be thought of as having more of an
elliptic character, and to a large extent have a perturbative
nature. The quadratic terms
\[
\P\left(  [A^k,  \p_j A_k] - 2 [ A^k, \partial_k A_j]\right) 
\]
can be thought of as the leading part of the nonlinearity. It is
crucial that these terms satisfy the cancellation property known as
the null condition.

As mentioned before, the Coulomb gauge works well for the small data problem (Theorem~\ref{t:coulomb}).
Concerning large data, however, one sees here that in order to properly set up
the Yang--Mills equation in the Coulomb gauge one would need to be
able to invert the operator $\partial^j \covD_j$.  Exactly the same
operator arises when one considers the linearization of the Coulomb
gauge condition. This works well in the small data problem, but not so
well for the large data problem.

\subsection{Local and global theory for the Yang--Mills heat flow} \label{subsec:caloric}
 Neglecting for the
moment the time component of the connection $A$, at fixed time we
consider the energy functional
\[
\spE(A_x) = \frac{1}{2} \int_{\R^4} \brk{F_{ij}, F^{ij}} dx.
\]
The \emph{Yang--Mills heat flow} is the gradient flow associated to this functional,
which has the expression 
\begin{equation}\label{caloric}
\partial_s A_i = \covD^{\ell} F_{\ell i}, \qquad A_i(s=0) = a_i.
\end{equation}
As written this system is invariant with respect to purely spatial
gauge transforms. To better frame the discussion, we  observe
that one can add a heat time component to the connection $A$ and
rewrite the Yang--Mills heat flow in a fully covariant fashion as
\begin{equation} \label{eq:cymhf}
F_{si} = \covD^\ell F_{\ell i}.
\end{equation}
Then one can view the  Yang--Mills heat flow equations in \eqref{caloric} as the 
effect of a gauge choice 
\[
A_s= 0,
\]
 (which we call the \emph{local caloric
gauge}) applied to the fully covariant Yang--Mills heat flow.
This is akin to using the temporal gauge for the hyperbolic Yang--Mills equation.

We start with the basic result:
\begin{theorem} \label{t:ymhf-loc}
The problem \eqref{caloric} is locally well-posed for data $a \in \dot{H}^{1}$.
\end{theorem}
The assumption $a \in \dot{H}^{1}$ restricts $a$ (and thus the solution) to the topologically trivial class.
This is natural in view of our goal of constructing the caloric gauge, and also for the eventual application to the Threshold Theorem (Theorem~\ref{thm:threshold}).

In the study of \eqref{caloric}, a key role is played by the
$L^3_{s,x}$ norm of the curvature $F_{ij}$.  Precisely, the solution
to \eqref{caloric} can be continued and uniform covariant parabolic estimates
for the solution can be proved for as long as $\|F\|_{L^3}$ remains
finite. This motivates the following definition
for the \emph{caloric size} of a connection $a$:
\[
\hM(a) = \left\{ 
\begin{array}{cl}
\int_{R^+ \times \R^4} |F(s,x)|^3 ds dx & \text{if the solution to \eqref{caloric} is global}, \cr
\infty & \text{otherwise}.
\end{array}
\right.
\]
We note that this is a scaling- and gauge-invariant quantity. 

As described below, the caloric gauge is
defined only for connections $a$ for which $\hM(a)$ is finite. This is
an open subset of $\dot H^1$, as $\hM(a)$ has a locally Lipschitz
dependence on $a$ whenever finite.  Furthermore, for such $a$ we can
describe the behavior of its Yang--Mills heat flow at infinity as
follows:
\begin{theorem}[\cite{OTYM1}]\label{t:global-heat}
Let $a \in \dot H^1$ be a connection so that $\hM(a) < \infty$. Then 
the corresponding solution has the property that the limit
\[
\lim_{s \to \infty} A(s) = a_\infty
\]
exists in $\dot H^1$. Further, the limiting connection is flat,
$f_\infty = 0$.
\end{theorem}

The main technical difficulty with \eqref{caloric} is that it is only degenerate parabolic.  
Precisely, \eqref{caloric} can be formally viewed as a coupling of a strongly parabolic system
for $F$ (which we think of as the curl of $A$) and a transport equation for the divergence of $A$.

We note that there is an alternate gauge choice which circumvents this issue, namely the \emph{de Turck
gauge} 
\[
A_0 = \rd^j A_j,
\]
 where the Yang--Mills heat flow becomes strongly
parabolic and is easier to solve locally. In our formalism, the classical de Turck trick of compensating the degeneracy by a suitable $s$-dependent gauge transformation amounts to solving \eqref{eq:cymhf} in this gauge, hence the name.

Unfortunately, the transition from local to global is impossible in the de Turck gauge; in other words, Theorem~\ref{t:global-heat} is false in the de Turck gauge. One can see this by considering the evolution of flat connections. This is
trivial under the local caloric gauge, but yields a $4+1$ dimensional
harmonic heat flow for maps into $\G$ in the de Turck gauge, which is known to possibly blow up.

Our approach is instead based on a version of the de Turck trick for the \emph{linearization} of \eqref{caloric} (namely, \eqref{caloric-lin} below). In this scheme, an auxiliary flow called the \emph{dynamic Yang--Mills heat flow} plays a major role. We will return to discussion of this idea in Section~\ref{subsec:dymhf}.

For now, we proceed to describe our next result proved in \cite{OTYM1}, which asserts that all connections with 
energy below threshold $2\Egs$ have finite caloric size, and thus Theorem~\ref{t:global-heat} applies:

\begin{theorem}[Threshold Theorem for the heat flow] \label{t:caloric}
There exists a non\-decreasing function
\[
\hM: [0,2\Egs) \to [0,\infty)
\]
so that for every connection $1$-form $a \in \dot H^1$ with subthreshold energy  $\nE < 2 \Egs$, 
we have 
\begin{equation}
\hM(a) \leq \hM(\nE)
\end{equation}
\end{theorem}

 This is proved using a concentration compactness type
argument. The key ingredient is the \emph{energy monotonicity formula}
\begin{equation*}
\spE(A(s_{1})) - \spE(A(s_{2})) = - \int_{s_{1}}^{s_{2}} \int \brk{\covD^{\ell} F_{\ell j}, \covD^{k} \tensor{F}{_{k}^{j}} } \, dx ds.
\end{equation*}
This formula yields good control of $A$ in the local caloric gauge, but not in the de Turck gauge.
The same argument also gives the corresponding Dichotomy Theorem:

\begin{theorem} [Dichotomy Theorem for the heat flow]\label{thm:dichotomy-heat}
For any $a \in \dot{H}^{1}$, one of the following two properties
  must hold for the maximally extended solution:
\begin{enumerate}[label=(\roman*)] 
\item The solution is topologically trivial, global and $\hM(a) < \infty$. 

\item The solution bubbles off a harmonic Yang--Mills connection either
\begin{enumerate}[label=(\alph*)] 
\item at a finite blow-up time, or

\item at infinity.
\end{enumerate}
\end{enumerate}

\end{theorem}

The bubbling argument here has roots in the classical work of Struwe \cite{Str} (see also Schlatter \cite{Sch1}) on compact manifolds.
In comparison, the significance of the above theorems lies in the precise asymptotics of the Yang--Mills heat flow on the noncompact space $\bbR^{4}$, which allows us to construct the caloric gauge.

\subsection{Caloric connections and the caloric manifold}
Since the limiting connection $a_\infty$ given by Theorem~\ref{t:global-heat}
is flat, it must be gauge equivalent to the zero connection.
Precisely, there exists a gauge transformation $O$ with the property that
\[
a_{\infty,j} =O^{-1}  \partial_j O  .
\]
Here $O= O(a) \in \dot H^2$ (interpreted in the sense that $O_{;j} := \partial_j O  O^{-1} \in
\dot H^1$) is unique up to constant gauge transformations.  Conjugating the
full heat flow with respect to such an $O$ yields a gauge equivalent
connection
\[
\tA_j = O A_j O^{-1} - O_{;j} ,
\]
which solves the Yang--Mills heat flow, and satisfies $\ta_\infty = 0$.
This lead us to the following definition of \emph{caloric} connections:

\begin{definition} \label{d:C}
We will say that a connection $a \in \dot H^1$ is  caloric  if $a_\infty = 0$. 
We  denote the set of all such connections by $\calC$.
\end{definition}

Theorem~\ref{t:caloric} can then be restated as an existence result for gauge equivalent caloric connections:

\begin{theorem}[\cite{OTYM1}] \label{t:caloric-g}
  For every connection $a \in \dot H^1$ with $\hM(a) < \infty$ there
  exists a gauge equivalent caloric connection $\tilde a \in \dot
  H^1$, which is unique up to constant gauge transformations. In particular,
  this conclusion holds for all subthreshold connections.
\end{theorem}

The connection $\ta$ is defined as 
\[
\ta_j = O a_j O^{-1} - O_{;j}, \qquad O = O(a) .
\]
We note that the two connections have the same caloric size, $\hM(a)= \hM(\ta)$. 

To solve the Yang--Mills equation in the caloric gauge we need to view
the family $\calC$ of the caloric gauge connections with energy below the
ground state energy as an infinite dimensional manifold. Here the $\dot H^1$ 
topology is no longer sufficient, so we introduce the slightly stronger topology
\[
\bfH = \{ a \in \dot H^1: \ \partial^j a_j \in \ell^1 L^2\}
\]
which reflects the fact, to be discussed later in more detail, that caloric connections
satisfy a generalized,  nonlinear form of the Coulomb gauge condition.  Then we have
\begin{theorem} [\cite{OTYM1}] \label{t:smooth}
For any caloric subthreshold connections $a$ with energy $\nE$ and caloric size $\hM$ 
we have the $\bfH$ bound
\begin{equation}
\| a \|_{\bfH} \lesssim _{\nE,\hM} 1
\end{equation}
 The set $\calC$ of all $\dot H^1$ caloric connections is a $C^1$
 infinite dimensional submanifold of $\bfH$.
\end{theorem}

We denote 
\[
\ta = \Cal(a).
\]
For arbitrary subthreshold $a \in \dot H^1$ this is only defined as an
equivalence class, modulo constant conjugations. However, if in
addition we know that $a \in \bfH$, then $O(a)$ is continuous, and we
can fix its choice by imposing the additional condition
\begin{equation}\label{O-inf}
\lim_{x \to \infty} O(x) = Id.
\end{equation}
With this choice we have the following regularity property:

\begin{theorem}
  The map $a \to O(a)$ is continuous (though not
    Lipschitz) from $\dot H^1$ to\footnote{Here $\dot H^2$ needs to
    be interpreted as a quotient space, modulo constant conjugations}
  $\dot H^2$.  It is also locally $C^1$ from $\bfH$ to\footnote{Here
    the action of the group of constant conjugations can be eliminated
    by using the condition \eqref{O-inf}.} $\dot H^2 \cap C^0$.
\end{theorem}

\subsection{The tangent space and caloric data sets}
\label{subsec:caloric-pf}

Finite energy caloric Yang--Mills waves will be continuous functions of time
which take values into $\calC$.  They are however not smooth in time,
instead their time derivative will merely belong to $L^2$.  Because of
this, we need to take the closure of its tangent space $T\calC$ (which
a-priori is a closed subspace of $\bfH$)  in $L^2$.  This is denoted 
by $T_a^{L^2}\calC$. It is also convenient to have a direct way of 
characterizing this space; that is naturally done via the linearization of the caloric flow:

  \begin{definition}\label{d:tangent}
For a caloric gauge connection $a \in \calC$, we say that $ L^2  \ni b \in T_a^{L^2} \calC$
iff the solution to the linearized local caloric gauge Yang--Mills heat flow equation
\begin{equation} \label{caloric-lin}
\partial_s B_k = [B^j,F_{kj}] + \covD^j ( \covD_k B_j - \covD_j B_k), \qquad B_k(0) = b_k
\end{equation}
satisfies 
\[
\lim_{s \to \infty} B(s) = 0.
\]
\end{definition}

Turning our attention now to the Yang--Mills flow, we will now consider
solutions which at any fixed time $t$ are in the caloric gauge,
$A_x(t) \in C$. 
 \begin{definition} 
An initial data for the Yang--Mills equation in the caloric gauge is a pair 
$(a,b)$ where $a \in \calC$ and $b_k \in T_a^{L^2} \calC$. 
\end{definition}

The transition from one time to another requires
understanding the linearization of the Yang--Mills heat flow.
As in the Coulomb gauge, we will consider the spatial component of the 
connection as the dynamic variable, and view the temporal part of the connection
as an auxiliary variable. We begin our discussion by considering the initial data. 
To connect a  general initial data $(a_k,e_k)$ with caloric initial
data we have the following result:
\begin{theorem}\label{t:data}
\begin{enumerate}
\item For any initial data pair $(a, e) \in \dot H^1 \times L^{2}$ with finite caloric size, there exists a caloric gauge
    data set $(\tilde a, b) \in T^{L^{2}} \calC$ and $a_0 \in
    \dot H^1$, unique up to constant gauge transformations and with continuous dependence in this quotient topology, so that $(\tilde a, \tilde e)$ is
    gauge equivalent to $(a,e)$ and
    \[
    \tilde e_k = b_k - (\covD_{\tilde{a}})_{k} a_0 .
    \]
\item For any caloric gauge initial data set $(\tilde a, b) \in T^{L^{2}}
    \calC$, there exists a unique $a_0 \in \dot H^1$, with Lipschitz
    dependence on $(a, b) \in \dot H^1 \times L^2$, so that
    \[
    e_k = b_k - (\covD_{a})_k a_0
    \]
    satisfies the constraint equation
    \eqref{eq:YMconstraint}. 
\end{enumerate}
\end{theorem}

In view of this result, we can fully describe caloric Yang--Mills waves as continuous functions
\[
I \ni t \to (A_x(t), \partial_0 A_x(t)) \in T^{L^2} \calC.
\]
An important role in the proof of this theorem is played by the
following nonlinear div-curl type decomposition for the tangent space
$T_a^{L^2} \calC$:

\begin{theorem}  \label{thm:div-curl}
Let $a \in \calC$ with energy $\nE$ and caloric size $\hM$. Then for each 
$e \in L^2$ there exists a unique decomposition
\begin{equation}
e = b - \covD a_0, \qquad b \in T_a^{L^2} \calC, \qquad a_0 \in \dot H^1.
\end{equation}
with the corresponding bound
\begin{equation}
\| b\|_{L^2} +\|a_0\|_{\dot H^1} \lesssim_{\nE,\hM} \|e \|_{L^2}.
\end{equation}
\end{theorem}
 
Proving the latter theorem, in turn, requires understanding of the linearized equation \eqref{caloric-lin}; we will return to this issue in Section~\ref{subsec:dymhf}.
 
\subsection{The dynamic Yang--Mills heat flow and the hyperbolic Yang--Mills equation}
To proceed further, given a caloric Yang--Mills wave on $I$, we seek to interpret the (covariant) hyperbolic Yang--Mills equation
\begin{equation}\label{ym-cov}
\covD^\alpha F_{\alpha \beta} = 0, 
\end{equation}
as gauge dependent hyperbolic evolutions for $A_x$. 
Separating these equations into 
\begin{equation}\label{ym-x}
\covD^\alpha \covD_\alpha A_k = \covD^k \covD^\alpha A_\alpha - [A^\alpha, \covD_k A_\alpha],
\end{equation}
respectively,
\begin{equation}
\covD^k \covD_k A_0 = \covD_0 \covD^k A_k  - [A^k, \covD_0 A_k],
\end{equation}
we seek to interpret the first equation as a hyperbolic evolution for $A_x$,
and the second as an elliptic compatibility condition for $A_0$. This is achieved 
in several steps as follows:

\bigskip

(i) First, we show that the pair $(A_{x}, \rd_{0} A_{x}) \in T^{L^2} \calC$ 
satisfies a generalized Coulomb like condition,
\begin{equation} \label{eq:da}
\partial^k A_k = \DA(A), \qquad  \partial^k A_k = \DB(A,B),
\end{equation}
where $\DA$ and $\DB$ are nice maps on $T^{L^{2}} \calC$, which contains an explicitly computed quadratic part, as well as purely perturbative higher order terms. Of course, this step does not have to anything to do with \eqref{ym-cov}, and holds for any pair in $T^{L^{2}} \calC$.
The key computation for $\rd^{k} A_{k}$ is
\begin{equation*} 
	\rd^{k} A_{k} = - \int_{0}^{\infty} \rd^{k} \rd_{s} A_{k}(s) \, d s = - \int_{0}^{\infty} \covD^{k} F_{sk}(s) + \hbox{(quadratic and higher)}
\end{equation*}
but by \eqref{caloric}, the linear term vanishes. A similar computation holds for $\rd^{k} B_{k}$.

\bigskip

(ii) Next, we use the $\beta = 0$ part of the equation \eqref{ym-cov} to show 
that $A_0$ is uniquely determined by $A_x$ and $B_x = \partial_0 A_x$, i.e.,
\[
A_0 = \bfA_0(A_x,B_x)
\]
where $\bfA_0$ is a nice smooth map on $T^{L^2} \calC$ which contains
an explicitly computed quadratic part, as well as purely perturbative
higher order terms.

 \bigskip

(iii)  Moreover, we use the $\beta \neq 0$ part of the equation \eqref{ym-cov} to show 
that $\covD^0 A_0$ is uniquely determined by $A_x$ and $B_x = \partial_0 A_x$,
\[
\covD^0 A_0 = \DA_0(A_x,B_x)
\]
where $\DA_0$ is a nice smooth map on $T^{L^2} \calC$ which again contains
an explicitly computed quadratic part, as well as purely perturbative
higher order terms.

\bigskip

The above steps allow us, just as in the case of the Coulomb gauge, to view the spatial 
part of the connection $(A_x,\partial_0 A_x) \in T^{L^2} \calC$ as the dynamical variable,
and $A_0$, $\partial_0 A_0$ as dependent variables.  
Precisely, we can recast the equations \eqref{ym-x} in the form
\begin{equation}\label{ymcg-main}
\Box_A A_k = \P [A_{x},\partial_k A_{x}] 
+  2\Delta^{-1} \partial_k \bfQ(\partial^\alpha A_{x},\partial_\alpha A_{x}) 
+ R(A,\partial_t A),
\end{equation}
where $[A_{x}, B_{x}]$ is a shorthand for $[A^{\ell}, B_{\ell}]$, 
and $\bfQ$ is a symmetric bilinear form with symbol\footnote{
Given a scalar-valued symbol $m(\xi, \eta)$, our definition of the associated bilinear multiplier is
\begin{equation*}
	\iint e^{i x \cdot (\xi + \eta)} m(\xi, \eta) [\hat{A}_{x}(\xi), \hat{B}_{x}(\eta)] \, \frac{d \xi}{(2 \pi)^{4}} \frac{d \eta}{(2 \pi)^{4}}.
\end{equation*}
} 
\[
\bfQ(\xi,\eta) = \frac{\xi^2-\eta^2}{2(\xi^2+\eta^2)}.
\]
Here on the right we have two explicit quadratic terms depending only on $A_x$,
and its time derivative, both of which have a favorable null structure, and 
a remainder higher order term $R$ which admits favorable $L^1L^2$ bounds
and thus only plays a perturbative role. However, in the covariant d'Alembertian
on the left, we still have the coefficients $A_0$ and $\covD_0 A_0$, which are determined 
as above in terms of  $A_x$ and $\partial_t A_x$:
\begin{equation} \label{a0-quad-}
\begin{split}
A_0 = & \  \bfA_0(A_x,B_x) =        \bfA_0^2(A_x,B_x) +  \bfA_0^3(A_x,B_x), \\
\covD_0 A_0 = & \ \DA_0(A_x,B_x) =   \DA_0^2(B_x,B_x) +  \DA_0^3(A_x,B_x),
\end{split}
\end{equation}
Here the quadratic terms  $\bfA_0^2(A_x,B_x)$, $\DA_0^2(A_x,B_x) $ 
are explicit translation invariant bilinear  forms,
\begin{align} 
\bfA_0^2(A_x,B_x) =& \Delta^{-1} [A_x,B_x] + 2 \Delta^{-1} \bfQ(A_x,B_x), \label{eq:a0-quad} \\ 
\DA_0^2(B_x,B_x)=& - 2 \Delta^{-1} \bfQ(B_x,B_x). \label{eq:da0-quad}
\end{align}
   The remainders
 $\bfA_0^3(A_x,B_x)$, $ \DA_0^3(A_x,B_x)$, however, are not explicit but satisfy
favorable bounds. Of these only the  quadratic part of $A_0$ plays a nonperturbative role.

Finally, $A_{x}$ is also subject to a compatibility condition
\begin{equation}\label{da}
\partial^k A_k = \DA(A) :=  \bfQ(A, A) + \DA^{3}(A),
\end{equation}
where $\DA^{3}$ is perturbative.

To study the small data problem it would be sufficient to work with
the equation \eqref{ymcg-main}. However, for the large data problem we
also need to flow the wave equation in the parabolic direction, which in turn requires us to specify the $s$-evolution equation for $A_{0}$. Our choice is to use the \emph{dynamic Yang--Mills heat flow}
\begin{equation} \label{eq:dymhf}
	F_{s \alp} = \covD^{\ell} F_{\ell \alp},
\end{equation}
which is the (covariant) Yang--Mills heat flow \eqref{eq:cymhf} adjoined with $F_{s0} = \covD^{\ell} F_{\ell 0}$.

For nonzero heat-times $s$, \eqref{ym-cov} now becomes
\begin{equation}\label{ym-cov-w}
\covD^\alpha F_{\alpha \beta}(s) = w_\alpha,
\end{equation}
where in general $w_\alpha$, called the \emph{Yang--Mills tension field}, is nontrivial as the two flows (wave and heat) do not commute.
Thus additional steps are needed:

\bigskip

(iv)  We compute parabolic evolutions for $w_\alpha$, showing that at time $t$ they
depend only on the data $A_x(t)$, $\partial_t A_x(t)$ and of course on $s$,
\[
w_\alpha = \bfw_\alpha(A_x(t),\partial_t A_x,s).
\]
Moreover, we separate $\bfw_\alpha$ into an explicit quadratic part and 
a higher order term,
\[
\bfw_\alpha(s) = \bfw_\alpha^2(s) + \bfw_\alpha^3(s)
\]
where the latter is purely perturbative.

\bigskip

(v) Finally, we recalculate $A_0$ and $\covD^0 A_0$ to include the dependence on 
$w(s)$, and write the analogue of the equation \eqref{ymcg-main}
for $A_x(s)$,
\begin{equation}\label{ymcg-main-s}
\begin{split}
\Box_{A(s)} A_k(s) = & \  \P [A^j(s),\partial_k A_j(s)] 
+  2\Delta^{-1} \partial_k \bfQ(\partial^\alpha A^j(s),\partial_\alpha A_j(s))  
+ R(A(s),\partial_t A(s)) \\
& \ + \P \bfw_k^2(s) + R_s(A,\partial_t A)
\end{split}
\end{equation}
The extra terms on the right are matched by a like contribution to the quadratic part of $A_0$,
i.e. \eqref{a0-quad-} is replaced by
\begin{equation} \label{a0-quad}
A_0(s) = \bfA_0^2(A(s),B(s)) +   \bfA_0^3(A(s),B(s)) +  \Delta^{-1} \bfw_0^2(A,B) +  \bfA^3_{0;s}(A,B)
\end{equation}
The $s$ dependent terms in the above equations depend on the original connection $A$ and not just 
on $A(s)$. However, they have the redeeming feature that they are concentrated at a single dyadic frequency
$s^{-\frac12}$.

The analysis of the equation \eqref{ymcg-main-s} is now very similar
to that of \eqref{ymcg-main}, with the minor proviso that the
quadratic terms in $\bfw$ in the two equations above have a very mild
nonperturbative role, and exhibit a null form type cancellation.

\subsection{Remarks on the dynamic Yang--Mills heat flow} \label{subsec:dymhf}
 In \cite{OTYM1}, the dynamic Yang--Mills heat flow \eqref{eq:dymhf} plays a major role in our proofs in several different ways:
\begin{enumerate}[label=(\roman*)]
\item {\it As a gauge covariant smoothing flow for spacetime
    connections.}  This is the most direct interpretation; \eqref{eq:dymhf} was used in this capacity
to fix the evolution of $w_{\mu}(s)$ in the preceding subsection. 

\item {\it As a means to perform the ``infinitesimal de Turck trick''
    for the linearized Yang--Mills heat flow in the local caloric
    gauge.}  As alluded to earlier, our understanding 
    of \eqref{caloric} is based on its linearization \eqref{caloric-lin},
    which in turn is analyzed through a version of de Turck trick. It implemented as follows, using \eqref{eq:dymhf} as a useful auxiliary tool:
    \begin{itemize}
\item   Given a one-parameter family of Yang--Mills heat flows
  $A_{j}(t, x, s)$ with data $a_{j}(t, x)$ $(t \in I, x \in \bbR^{4},
  s \in J)$, we add a $t$-component $A_{0}(t, x, s)$ and view it as a
  connection 1-form on $I \times \bbR^{4} \times J$.
  In the $s$-direction, we then \emph{impose} the dynamic Yang--Mills
  heat flow \eqref{eq:dymhf}.

\item  Then the key idea is to work with
  \begin{equation} \label{eq:F0j-deT} F_{0j} = \rd_{t} A_{j} -
    \covD_{j} A_{0}.
  \end{equation}
  As opposed to $\rd_{t} A_{j}$, which solves \eqref{caloric-lin},
  $F_{0j}$ has the advantage of obeying a \emph{nondegenerate}
  covariant parabolic equation:
  \begin{equation*}
    \covD_{s} F_{0j} - \lap_{A} F_{0j} - 2 ad(\tensor{F}{_{j}^{\ell}})F_{0\ell} = 0.
  \end{equation*}
  Solving this equation would determine $F_{0j}$ from any data
  $F_{0j}(s=0) = e_{j}$.  We \emph{choose} $e_{j} = \rd_{t} a_{j}$,
  which amounts to prescribing $a_{0} = 0$.  Then $A_{0}$ may be
  determined by integrating $\rd_{s} A_{0} = F_{s0} =
    \covD^{\ell} F_{\ell 0}$, and then we come
  back to the solution $\rd_{t} A$ of \eqref{caloric-lin}.

\end{itemize}

\item {\it As a means to obtain useful representation of projection to
    the caloric manifold.}  This is a variant of (2).  Previously, we chose to initialize $a_{0} =
  0$. When $a(t=0)$ is a caloric connection, another natural choice is
  to set $A_{0}(s = \infty) = 0$, which amounts to requiring that the nearby $a(t)$'s are also
  caloric. Integrating $\rd_{s} A_{0} = \covD^{\ell} F_{\ell 0}$ from $s = \infty$ to $0$, we
  obtain 
  \begin{equation} \label{eq:a0-cal}
    a_{0} = - \int_{0}^{\infty} \covD^{\ell} F_{\ell 0}(s) \, ds.
  \end{equation}
  By \eqref{eq:F0j-deT}, we have
  \begin{equation*}
    e_{j} = \rd_{t} a_{j} - \covD_{j} a_{0}.
  \end{equation*}
  Since $a(t)$'s are caloric, $\rd_{t} a_{j}$ clearly belongs to
  $T_{a} \calC$, whereas $\covD a_{0}$
  is a pure covariant gradient. This procedure proves Theorem~\ref{thm:div-curl}, while yielding a useful representation formula \eqref{eq:a0-cal}.
\end{enumerate}

\section{Energy dispersed caloric Yang--Mills waves} \label{sec:ed}

Our second article \cite{OTYM2} is concerned with the hyperbolic
Yang--Mills equation in the caloric gauge, namely the equation
\eqref{ymcg-main} with the auxiliary variables $A_0$ and $D_0A_0$ as
in \eqref{a0-quad-} and the constraints \eqref{da}.

\subsection{Main results in the caloric gauge}
The first result is a local well-posedness result which uses the
notion of \emph{$\eps$-energy concentration scale}, defined as
\begin{equation*}
	r_{c}^\eps[a, e] = \sup \set{r : \sup_{x} \int_{B_{r}(x)} \abs{f}^{2} + \abs{e}^{2} \, \ud x \leq \eps^2}.
\end{equation*}
Then we have

\begin{theorem}[\cite{OTYM2}] \label{thm:lwp-simple} There exists a positive
  non-increasing function $\eps_{\ast}(\nE,\hM)$ so that for any
  initial data set $(a, e)$ with energy $\nE$ and initial caloric size
  $\hM$,  that the Yang--Mills equation in the caloric gauge is locally
well-posed in $\dot H^1 \times L^2$  on the time interval $[-r^{\eps_\ast}_{c}, r_{c}^{\eps_\ast}]$.
\end{theorem}
We omit here the precise meaning of well-posedness, and instead refer
the reader to Theorem~\ref{t:main} in the last section.  Precisely,
the conclusions of Theorem~\ref{t:main} hold restricted to the
interval $[-r_{c}^{\eps_{\ast}}, r_{c}^{\eps_{\ast}}]$.

The second main result in \cite{OTYM2} uses the notion of energy
dispersion, first introduced in \cite{ST1} in the Wave Maps
context. For a connection $A$ on a time interval $I$, we define its \emph{energy
  dispersion} as
\begin{equation*}
	\nrm{F}_{ED[I]}  = \sup_{k} 2^{-2k} \nrm{P_{k} F}_{L^{\infty} L^{\infty}[I]}.
\end{equation*}
Then we have:

\begin{theorem} \label{thm:ed-simple} There exists a positive
  non-increasing function $\eps(\nE)$ and a nondecreasing function
  $M(\nE)$ such that if $A$ is a caloric Yang--Mills wave on $I$ with
  energy $\nE$ and initial caloric size $\hM \lesssim_{\nE} 1$ so that
  $\nrm{F}_{ED} \leq \eps(\nE)$, then\footnote{The control norm $S^1$
    will be described shortly.}  $\| A\|_{S^1[I]} \leq M(\nE)$ and $A$
  can be continued (as a well-posed solution in the sense of
  Theorem~\ref{thm:lwp-simple}) past finite endpoints of $I$.
\end{theorem}

We also note that the initial assumption on $\hM$ only serves to
prevent it from being very large. With this assumption, we actually
show that $\hM(A) \ll 1$ in the entire interval $I$. By
Theorem~\ref{t:caloric}, this assumption can be entirely omitted for
subthreshold energies.

These theorems, or rather their contrapositives, can be considered as
continuation criteria for the hyperbolic Yang--Mills equation in the
caloric gauge. By providing an accurate description of how singularities
may occur, they furnish a starting point for the bubble extraction
argument in \cite{OTYM3}, as it will be explained in
Section~\ref{subsec:no-bubble}.

One downside of using either the caloric gauge (or the Coulomb gauge) is that causality is lost. 
To remedy this, we prove that the well-posedness property can be transferred from the caloric gauge to the temporal gauge
$A_0 = 0$. As a result, we obtain:

\begin{theorem} \label{thm:temp-small}
The hyperbolic Yang--Mills equation in $\bbR^{4+1}$   is globally well-posed in the temporal gauge for 
all initial data with small energy.   
\end{theorem}
Unlike the caloric gauge results, however, a downside of Theorem~\ref{thm:temp-small} is that it does not
provide the $S^1$ regularity of solutions, or any other dispersive bounds.

In the remainder of this section, we will give an overview of ideas
in the proofs of Theorems~\ref{thm:lwp-simple},
\ref{thm:ed-simple} and \ref{thm:temp-small}.

\subsection{Function spaces}
To state the results more precisely, and also to discuss their proof,
it is necessary to outline the function spaces framework used in
\cite{OTYM2}, whose main components are the same as in \cite{KST,
  KT}. The core solution space, which we denote by $S^{1}[I]$, is a
Banach space of functions on $I \times \bbR^{4}$ with the property
that elements of $S^{1}[I]$ inherit estimates satisfied by free waves
in the energy class (i.e., $\Box u = 0$ with $(u, \rd_{t} u)(0) \in
\dot{H}^{1} \times L^{2})$, such as energy estimates, Strichartz
estimates, (null form) bilinear estimates etc. The corresponding
nonlinearity space, denoted by $N[I]$, is defined, on the one hand,
small enough to satisfy the inhomogeneous estimate
\begin{equation} \label{eq:en-est-free}
	\nrm{u}_{S^{1}[I]} \aleq \nrm{(u, \rd_{t} u)(0)}_{\dot{H}^{1} \times L^{2}} + \nrm{\Box u}_{N[I]},
\end{equation}
and on the other hand, large enough to contain (at least, most of) the
nonlinearities of the wave equation \eqref{ymcg-main}.

Construction of these spaces builds up on many prior works. 
The space $N[I]$ is simply the sum of
the dual energy space (i.e., $L^{1} L^{2}[I]$) and a dual $X^{s, b}$
space.
Building blocks of the space $S^{1}[I]$ include the energy space (i.e.,
$\nrm{\nb u}_{L^{\infty} L^{2}[I]}$), the Strichartz spaces (i.e.,
$\nrm{\abs{D}^{-\alp} \nb u}_{L^{p} L^{q}[I]}$ with admissible $\alp,
p, q$), an $X^{s, b}$ space \cite{KlMa0, Bour}, the refined Strichartz
spaces with radial frequency localization \cite{Kl-Tat}, and the null
frame space \cite{Tat, Tao2}.  Moreover, we also add a new component 
$S^{sq}$ (to be described in Section~\ref{subsec:temp-small}), which is used in the proof of Theorem~\ref{thm:temp-small}.
 For the precise definition, we refer to
\cite[Section~4]{OTYM2}.

The $S^{1}[I]$-norm serves the role of a controlling norm for the
caloric Yang--Mills waves. More precisely, we show in \cite{OTYM2}
that finiteness of this norm implies finer properties of the solution
itself and those nearby, such as frequency envelope control,
persistence of regularity and scattering for $A_{x}$, as well as weak
Lipschitz dependence and local-in-time continuous dependence for the
nearby solutions. For details, see the structure theorems in
\cite[Section~4]{OTYM2}.

\subsection{Truncated energy dispersion and the central result}
It turns out that Theorems~\ref{thm:lwp-simple} and
\ref{thm:ed-simple} can be proved essentially at the same time. The
idea is to use smallness of the \emph{truncated energy dispersion} at
frequencies higher than $2^{m}$,
\begin{equation}\label{eq:ed-m}
\nrm{F}_{ED_{>m}[I]} = \sup_{k > m}  2^{-2k} \nrm{P_{k} F}_{L^{\infty}L^{\infty}[I]},
\end{equation}
matched with shortness of the time interval on the scale $2^{-m}$.
The central result of \cite{OTYM2} reads as follows.
\begin{theorem} \label{thm:ed-main} There exist a non-decreasing
  positive function $M(\nE,\hM)$ and non-increasing positive functions
  $\eps(\nE,\hM)$ and $T(\nE,\hM)$, so that the following holds: For
  all regular subthreshold caloric Yang--Mills waves $A$ in a time
  interval $I$ with energy $\nE$ and initial caloric size $\hM$, if we
  have
\begin{equation} \label{eq:ed-hyp}
 \quad \nrm{F}_{ED_{\geq m}[I]} \leq \eps(\nE,\hM), \quad
  \abs{I} \leq 2^{-m} T(\nE,\hM),
\end{equation}
then we must also have
\begin{equation} \label{eq:ed-bnd}
\nrm{A}_{S^{1}[I]} \leq M(\nE,\hM).
\end{equation} 
\end{theorem}

On the one hand, this theorem implies an $S^{1}[I]$ control norm bound
on a time interval of size $\leq 2^{-m}$ for data with sufficiently
small energy at frequencies $> 2^{m}$ (i.e., $\nrm{P_{>m} (A_{x},
  \rd_{t} A_{x})(0)}_{\dot{H}^{1} \times L^{2}}$ is small), which is
the case for data with energy concentration scale $\ageq 2^{-m}$. On
the other hand, it also implies an $S^{1}[I]$-bound, independent of
$I$, if the solution has small untruncated energy dispersion
$\nrm{F}_{ED[I]}$. As discussed above, these $S^{1}[I]$-norm bounds
prove Theorems~\ref{thm:lwp-simple} and \ref{thm:ed-simple},
respectively.

\subsection{Review of the small energy case: Perturbative 
nonlinearities and parametrix construction}
We begin with a brief discussion of the small energy case, where the
goal is to prove $\nrm{A_{x}}_{S^{1}[\bbR]}^{2} \leq C \nE$ for
sufficiently small $\nE$. This was carried out in \cite{KT}, which can
be viewed as one of the predecessors to this work, in the closely
related context of the Coulomb gauge\footnote{While the analysis in
  \cite{KT} is carried out in the Coulomb gauge
  $\rd^{\ell} A_{\ell} = 0$, it is not very different in the caloric
  gauge, as this also satisfies some form of generalized Coulomb
  condition $\partial^{\ell} A_{\ell} = \DA(A)$.  }.

The first step was to try to view the wave equation for $A_{x}$ as a perturbation of
the constant coefficient wave equation $\Box A_{x} = 0$.  While this is not possible,
we can view most of the nonlinearity as perturbative, and
estimate them in the space $N$. In this process, the primary (bilinear) null structure of the Yang--Mills equation, uncovered in \cite{KlMa2}, plays an essential role.
This leaves us with a single nonperturbative term, which arises in a
paradifferential fashion,
\begin{equation}\label{para}
(\Box + \Diff^{0}_{\P A}) A_{x} := \left( \Box + 2 \sum_{k} ad(P_{<k} \P^{\alp} A) \partial_\alpha  P_{k}\right)  A_{x} =  G
\end{equation}
where $\P_{x} A$ is the Leray projection of $A_{x}$, $\P_{0} A = A_{0}$ and $ G$ represents a nonlinear but perturbative contribution (which is small thanks to smallness of energy). 

Then the key step in \cite{KT} was to construct a parametrix for the
paradifferential operator $\Box + \Diff^{0}_{\P A}$, and prove that this parametrix
satisfies a good $N \to S^1$ bound akin to \eqref{eq:en-est-free}.
The rough idea is to try to find a gauge transform $O$ which renormalizes $\Box + \Diff^{0}_{\P A}$ to $\Box$ modulo a better behaved error, i.e., schematically
\begin{equation} \label{eq:box-p-err}
	(\Box + \Diff^{0}_{\P A}) Ad(O) - Ad(O) \Box= \hbox{(error)},
\end{equation}
and produce a parametrix by conjugating the constant coefficient solution operator by $Ad(O)^{-1}$.

This idea was indeed viable in the case of wave maps \cite{Tao2, ST1}, but not for Yang--Mills or Maxwell--Klein--Gordon (which may be regarded as a simpler model for Yang--Mills). The difference stems from the structure of the curvature $F[\P A]$, which is a geometric obstruction for gauge transformation of $A$ to $0$. Whereas the curvature depends at least quadratically on the solution in the case of wave maps, it is linear (to the leading order) in the solution $A$ for Yang--Mills or Maxwell--Klein--Gordon.

The way out of this difficulty was to consider instead an $Ad(\G)$-valued \emph{pseudodifferential} renormalization operator $Op(Ad(O))$. Heuristically, this generalization allows for separate renormalization of each plane wave solution, which is possible since it only oscillates in a single direction\footnote{This procedure eludes the geometric obstruction mentioned above, since curvature, being a $2$-form, always vanishes when restricted to a one-dimensional subspace. }. Using smallness of energy, it was shown that the parametrix obeys the desired $N \to S^{1}$, and also that the error in \eqref{eq:box-p-err} is perturbative. We remark that in the error estimate, not only the primary but also the secondary (trilinear) null structure, analogous to that in Maxwell--Klein--Gordon discovered in \cite{MS}, is crucial.

\subsection{Parametrix construction in the large energy case}
The difference in the large energy case is that we can no longer use
smallness of energy to control neither the perturbative part, nor the parametrix
for the paradifferential problem.  Thus, in order to be able to close
our estimates, we need to have new proxies for smallness.

We start with the paradifferential problem. In departure from the small energy case, but similar to \cite{ST1, OT2}, we introduce the \emph{large frequency gap} $\kpp \gg 1$ and consider the paradifferential operator
\begin{equation*}
	\Box + \Diff^{\kpp}_{\P A} = \Box + 2 \sum_{k} ad(P_{<k-\kpp} \P^{\alp} A) \rd_{\alp} P_{k},
\end{equation*}
where $A_{x}$ be a caloric Yang--Mills with finite $S^{1}[I]$-norm. The goal is to establish an $N \to S^{1}$ bound of the form
\begin{equation} \label{eq:paradiff-goal}
	\nrm{u}_{S^{1}[I]} \aleq_{\nrm{A_{x}}_{S^{1}[I]}} \nrm{(u, \rd_{t} u)(0)}_{\dot{H}^{1} \times L^{2}} + \nrm{(\Box + \Diff^{\kpp}_{\P A}) u}_{N[I]}.
\end{equation}

The proof proceeds by a parametrix construction, in a similar manner as \cite{KT}. However, the necessary smallness for proving the $N \to S^{1}$ bound for the parametrix now comes from taking the frequency gap $\kpp$ sufficiently large compared to $\nrm{A_{x}}_{S^{1}[I]}$. Moreover, to control the error, we rely on divisibility\footnote{That is, $I$ can be split into a controlled number of subintervals, on each of which the restricted norm is arbitrarily small.} of an appropriate weaker norm $\nrm{A_{x}}_{DS^{1}[I]}$ than $\nrm{A_{x}}_{S^{1}[I]}$.

\subsubsection*{Treating perturbative nonlinearity: Small energy dispersion and short time interval}
For the perturbative nonlinearity, smallness may be obtained via truncated energy dispersion and the length of $I$.
Roughly speaking, any unbalanced or close-angle frequency interaction is small (exponentially in the frequency ratio) for such nonlinearities, while balanced and far-angle interactions are controlled by $\nrm{F}_{ED_{>m}[I]}$ at frequencies $\ageq 2^{m}$, and by $2^{m} \abs{I}$ at frequencies $\aleq 2^{m}$.
In sum, we have
\begin{equation*}
\nrm{F}_{ED_{>m}[I]} \leq \veps, \quad 2^{m}\abs{I} \leq \veps \implies \nrm{(\Box + \Diff^{\kpp}_{\P A}) A_{x}}_{N[I]} \aleq_{\nrm{A_{x}}_{S^{1}[I]}} 2^{C \kpp} \veps^{\dlt}.
\end{equation*}

Unfortunately, this bound is insufficient for proving
Theorem~\ref{thm:ed-main}. The reason is that the $N \to S^{1}$ bound
\eqref{eq:paradiff-goal} for the paradifferential operator already
depends on the $S^{1}[I]$-norm of $A_{x}$, which is what we wish to
bound!

\subsection{Induction on energy} 
In order to break the circular argument, we perform an induction on
energy, following the scheme developed in \cite{ST1}. Roughly
speaking, the main idea is to view $A$ as a perturbation of another
solution $\tA$, which has a lower (linear) energy and hence obeys an
$S^{1}$-norm bound by an induction hypothesis. To make this idea work,
we need to carefully construct $\tA$ so that we may control the
difference $A - \tA$.

A preliminary step here is to show that  $\hM$ is essentially
conserved for solutions with small energy dispersion. Once this is
done, $Q$ becomes a fixed parameter and is omitted from the subsequent
discussion.

The induction argument is set up as follows, in terms of the linear
energy $E$ rather than the nonlinear one $\nE$.  The initial step is
provided by the small energy case, which proves \eqref{eq:ed-bnd} up
to sufficiently small $E > 0$, with $M(E) = C \sqrt{E}$ and any
choices of $\eps(E)$, $T(E)$. As the induction hypothesis, we assume
that there exist functions $\eps(\cdot)$, $T(\cdot)$ and $M(\cdot)$
such that \eqref{eq:ed-bnd} holds up to some $E$. Then the goal is to
extend these functions so that \eqref{eq:ed-bnd} holds up to
$E + c_{0}$ for some $c_{0} = c_{0}(E) > 0$. An essential point for
continuing this induction argument (in order to cover all subthreshold
solutions) is to ensure that the increment $c_{0}(E)$ is
\emph{independent} of the functions $\eps(\cdot)$, $T(\cdot)$ and
$M(\cdot)$ given by the induction hypothesis\footnote{Meanwhile,
  $\eps = \eps(E + c_{0})$, $T = T(E + c_{0})$ and $M = M(E + c_{0})$
  may (and indeed do) depend on $\eps(E)$, $T(E)$ and $M(E)$. We are
  allowed to choose these parameters in the order
  $c_{0} \to M \to T, \eps$.}.

We define $\tA$ by first flowing the data $\tA_{x}(0)$ and $\rd_{t}
\tA_{x}(0)$ by the Yang--Mills heat flow and the linearized
Yang--Mills heat flow, respectively, for some heat-time $s_{\ast}$,
then solving the Yang--Mills equation in caloric gauge in time. Taking
$\eps$, $T$ and $c_{0}$ sufficiently small, and choosing $s_{\ast}$
appropriately, we aim for the following two goals:
\begin{enumerate}[label=\roman*)]
\item $\tA$ exists on $I$ and $\nrm{\tA}_{S^{1}[I]} \leq M(E)$;
\item $\nrm{A - \tA}_{S^{1}[I]} \aleq_{M(E)} 1$.
\end{enumerate}

The cutoff heat-time $s_{\ast}$ can be chosen so that either a)
$s_{\ast} \ll 2^{-m}$ and $\nrm{\nb \tA(0)}_{L^{2}} = E$, or b)
$s_{\ast} \simeq 2^{-m}$ and $\nrm{\nb \tA(0)}_{L^{2}} \geq E$. In
both cases, provided that $\eps, T$ are sufficiently small, it can be
shown that $\tA_{x}$ is close to the Yang--Mills heat flow
$A_{x}(s_{\ast})$ of $A_{x}$. In Case~a), taking $\eps$ smaller if
necessary, we may ensure that $\nrm{\tF}_{ED_{\geq m}} \leq \eps(E)$
and Goal~i) follows from the induction hypothesis. In Case~b),
$\tA(0)$ is sufficiently smooth so that the desired conclusion can be
proved simply by higher order local well-posedness.

To accomplish ii), we need several ideas. First, we observe that the
linear energies $\nrm{\nb A_{x}(t)}_{L^{2}}$, $\nrm{\nb
  \tA_{x}(t)}_{L^{2}}$ of the solutions $A, \tA$ are conserved in $t$,
up to an error that can be made arbitrarily small by taking $\eps, T$
small enough.  Moreover, since $\tA$ is close to $A(s_{\ast})$, which
in turn is (at least heuristically) a low frequency truncation of $A$,
the frequency supports of $A - \tA$ and $\tA$ are essentially
separated. Therefore, approximate conservation of linear energies for
$A$ and $\tA$ implies
\begin{equation} \label{eq:approx-en-diff} \sup_{t \in I} \nrm{\nb
    (A_{x} - \tA_{x})(t)}_{L^{2}} \aleq_{E} \nrm{\nb A_{x}(0)}_{L^{2}}
  - \nrm{\nb \tA_{x}(0)}_{L^{2}} \leq c_{0}.
\end{equation}
To upgrade this to an $S^{1}[I]$-norm bound, we establish \emph{weak
  divisibility} of the $S^{1}$-norm of $\tA$, i.e., that we can split
$I = \cup_{k = 1}^{K} I_{k}$ so that
\begin{equation} \label{eq:weak-div} \nrm{\tA_{x}}_{S^{1}[I_{k}]}
  \aleq_{E} 1, \quad K \aleq_{M(E)} 1.
\end{equation}
Now viewing $A = \tA + (A - \tA)$ as a perturbation of $\tA$ on each
$I_{k}$, where the data for $A - \tA$ are reinitialized on each
interval using \eqref{eq:approx-en-diff}, we may bound the
$S^{1}$-norm of $A - \tA$ on each $I_{k}$ provided that $c_{0}$ is
small enough compared to the implicit constants in
\eqref{eq:approx-en-diff} and \eqref{eq:weak-div}. Importantly, these
are independent of $M(E)$! Thus Goal B follows by summing up these
bounds in $k = 1, \ldots, K$.

\subsection{Passing to the temporal gauge} \label{subsec:temp-small}
Finally, we describe the ideas behind the proof Theorem~\ref{thm:temp-small}. 
We wish to estimate the gauge tranformation $O$ from the caloric gauge into the temporal gauge, which solves the nonlinear transport equation
\begin{equation*}
	O^{-1} \rd_{t} O = A_{0}.
\end{equation*}
For $O$ to preserve $\dot{H}^{1}$ regularity of $A_{x}$, we need:
\begin{equation} \label{eq:temp-goal}
	\lap A_{0} \in \ell^{1} L^{2}_{x} L^{1}_{t}.
\end{equation}

The proof of \eqref{eq:temp-goal} relies on two observations. 
\begin{enumerate}[label=(\roman*)]
\item We note that the following \emph{square function norm} can be added to the $S^{1}$ norm, i.e.,
\begin{equation*}
	\nrm{\nb A_{x}}_{S^{sq}} \aleq \nrm{A_{x}}_{S^{1}}.
\end{equation*}
where
\begin{equation*}
\nrm{u}_{S^{sq}} = \nrm{\abs{D}^{-\frac{3}{10}} u}_{\ell^{2} L^{\frac{10}{3}}_{x} L^{2}_{t}}.
\end{equation*}
The relevance of $p = \frac{3}{10}$ is that it is the dual Stein--Tomas exponent for Fourier restriction to $\bbS^{3} \subseteq \bbR^{4}$. Indeed, the (adjoint) Stein--Tomas restriction theorem and Plancherel in time leads to
\begin{equation*}
	\nrm{e^{\pm i t \abs{D}} u}_{S^{sq}} \aleq \nrm{u}_{L^{2}},
\end{equation*}
which implies $\nb u \in S^{sq}$ for $\dot{H}^{1}$ free waves.
We extend this estimate to our parametrix, which allows us to add $S^{sq}$ into our $S^{1}$ norm.

\item In an order zero bilinear expression of the form $\bfO(A_{x}, \rd_{t} A_{x})$, the worst case is when $\rd_{t} A_{x}$ has the higher frequency.
Indeed, the ordinary product $[A_{x}, \rd_{t} A_{x}]$ fails to belong to $\ell^{1} L^{2}_{x} L^{1}_{t}$ because of this interaction. 
However, from \eqref{eq:a0-quad}, we see that the symbol of $\lap\bfA_{0}^{2}$ is
\begin{equation*}
	\lap \bfA_{0}^{2}(\xi, \eta) = \frac{2 \abs{\xi}^{2}}{\abs{\xi}^{2} + \abs{\eta}^{2}}.
\end{equation*}
which exhibits a favorable gain in the problematic $low \times high$ interaction!
\end{enumerate}

\section{Large data, causality and the 
temporal gauge} \label{sec:topo}

Unlike the first two papers, the third one \cite{OTYM2.5} is concerned
with large data solutions which are not necessarily topologically
trivial, and thus cannot be directly studied using the global caloric
gauge. The goal of \cite{OTYM2.5} is two-fold:
\begin{itemize}
\item To describe finite energy initial data sets topologically and analytically.

\item To use the temporal gauge in order to provide a good local theory for finite energy solutions.
\end{itemize}

For simplicity we will  work in two settings:

\begin{enumerate}[label=\alph*)]
\item  For initial data in $\bbR^4$ and solutions in $\R^{4+1}$, or time sections thereof.

\item For initial data in a ball $B_R$ and solutions in the corresponding uniqueness cone
$\calD(B_R) = \{ |x|+|t| < R\}$
or time sections thereof.
\end{enumerate}
In terms of the initial data, in addition to the energy, a key role is played by 
the energy concentration scale\footnote{For a singlet $a$, we define $r^{\eps}_{c}$ and $R^{\eps}_{c}$ by taking $e = 0$.}
\[
r_c^{\eps} = \sup  \{ r > 0 : \calE_{B_{r}(x) \cap X}(a,e)] \leq \eps \ \hbox{ for all } x \in X \},
\]
where $X = B_{R}$ or $\bbR^{4}$, as well as the outer concentration radius
\[
R_c^{\eps} =  \inf  \{ r > 0: \calE_{B(x,r)}[(a,e)] \leq \eps  \text{ for some } x \in \bbR^4\  \}.
\]

\subsection{Initial data surgery} 
Here we discuss a technical tool introduced in \cite{OTYM2.5}, which may
be of independent interest.  At various points in the analysis, we
need to perform a physical space localization of the Yang--Mills
solution.  By finite speed of propagation, this task amounts to
smoothly cutoff an initial data set $(a, e)$.  which turns out to be
nontrivial due to the presence of the constraint equation
\eqref{eq:YMconstraint}.  To address this issue, we prove the
following result:
\begin{theorem} \label{thm:div-solve} Let $B = B_{R_{0}}(0)$ be a ball
  centered at $0$, and let $a$ be a $\dot{H}^{1}$ connection on
  $\bbR^{4} \setminus B$. Then there exists a solution operator $h
  \mapsto e = T_{a} h$ to the equation
\begin{equation} \label{eq:div-eq-extr}
	\covD^{\ell} e_{\ell} = h \quad \hbox{ in } \bbR^{4} \setminus B, 
\end{equation}
with the following properties:
\begin{enumerate}
\item Boundedness: The operator $T_{a}$ is bounded from $\dot{H}^{-1}$
  to $L^{2}$, with a norm depending only on $\nrm{a}_{L^{4}}$.
\item Higher regularity: If $a$ and $h$ are smooth, then $T_{a} h$ is also smooth.
\item Exterior support: For any $R \geq R_{0}$, if $h = 0$ in $B_{R}(0)$, then $T_{a} h = 0$ in $B_{R}(0)$.
\end{enumerate}
\end{theorem}

In the case $a = 0$, \eqref{eq:div-eq-extr} becomes the usual
divergence equation and a desired solution operator $T_{0}$ may be
constructed explicitly. Exploiting the exterior support property of
$T_{0}$, $T_{a}$ is constructed in an essentially inductive manner,
starting from an annulus around $B$ (where $a$ can be treated
perturbatively) and proceeding outward.

As a quick corollary of Theorem~\ref{thm:div-solve}, we obtain the following initial data excision result.
\begin{proposition} \label{prop:id-excise} Let $(a, e)$ be a small
  energy data set in $B_{4} \setminus B_{1}$. Then 
\begin{enumerate}
\item We can find a small energy exterior data set $(\tilde{a},
  \tilde{e})$ in $\bbR^{4} \setminus B_{1}$, which agrees with $(a,
  e)$ in $B_{2} \setminus B_{1}$. Furthermore, if $(a, e)$ is smooth
  then $(\tilde{a}, \tilde{e})$ is also smooth.
\item We can find a small energy exterior data set $(\tilde{a},
  \tilde{e})$ in $\bbR^{4} \setminus B_{1}$, which is gauge equivalent
  to $(a, e)$ in $B_{4} \setminus B_{2}$. Furthermore, if $(a, e)$ is
  smooth then $(\tilde{a}, \tilde{e})$ is also smooth. 
\end{enumerate}
\end{proposition} 
The idea of the proof is to first naively extend $(a, e)$ to $\bbR^{4}
\setminus B_{1}$. This generates an error in the
constraint equation, which can be removed by applying
Theorem~\ref{thm:div-solve}.

\begin{remark} 
  Theorem~\ref{thm:div-solve} can clearly be generalized to other
  regularities and dimensions.  In particular, the operator $T_{a} :
  \dot{H}^{-1}(\bbR^{3} \setminus B) \to L^{2}(\bbR^{3} \setminus B)$
  can be used to prove an excision result for finite energy data on
  $\bbR^{3}$. We note that this furnishes an alternative approach to
  constructing local Coulomb gauges \cite{KlMa2} that avoids the need
  to prescribe boundary values.
\end{remark}

\subsection{Good global gauges}
In view of the gauge independence property, having control of the
energy of a connection $a$ says little about the $\dot H^1 \cap L^4$
size of $a$. This issue can sometimes be addressed by choosing a good
gauge, such as the local Coulomb gauge in Uhlenbeck's Lemma for small
energies, or the caloric gauge for subthreshold energies, see
Theorems~\ref{t:caloric}, \ref{t:smooth}. However, what if our
connection has larger energy?

We begin our discussion with initial data sets in a ball. In addition
to the energy $\nE$, we also use a second parameter, namely the energy
concentration scale $r_c = r_c^\eps$, with a small universal constant
$\eps$.  Then we have:

\begin{proposition} \label{prop:goodrep-ball}
 Given a connection $a$ in $B_R$ with energy $\nE$ and energy concentration scale $r_C$, there exists 
a gauge equivalent connection $\ta$ in $B_R$ which satisfies the bound
\begin{equation}
\| \ta\|_{\dot H^1 \cap L^4} \lesssim_{E, \frac{r_c}{R}} 1
\end{equation}
\end{proposition}

%

Also for initial data in $\R^4$, we also can find a good global gauge:

\begin{theorem} [Good global gauge] \label{thm:goodrep} Let
 $a \in  H^1_{loc}(\R^4)$ be a finite energy connection.  Then
  there exists a gauge equivalent representation $\ta$ of $a$  such that
\begin{equation*}
	\ta = - \chi O_{(\infty); x} + b
\end{equation*}
where $O_{(\infty)} (x)$ is a smooth $0$-homogeneous map taking values
in $\G$ and $B \in \dot H^1$. 
\end{theorem}

Finally, we remark on the relationship between Theorem~\ref{thm:goodrep} and topological classes of finite energy connections.
Precisely, the topological class of a connection $a$ can be parametrized by the 
homotopy class $[O]$ of the map $O$ in the above theorem, viewed as a map 
\[
O: \mathbb S^3 \to \G.
\]

\subsection{The temporal gauge and causality}

While we are not able to carry out the full analysis for the Yang--Mills equation in the temporal gauge,
we are nevertheless making  good use of it in our papers in an auxiliary role. This is due to the
following three properties:

\begin{enumerate}[label=(\roman*)]
\item  Local well-posedness for regular data. 

\item Causality, i.e. finite speed of propagation

\item Agreement with caloric gauge at the linear level.
\end{enumerate}

In our sequence of papers we are taking advantage of these three
properties at different places in the analysis. Property (i), for
instance, is used in order to prove a local well-posedness for regular
data in the caloric gauge, simply by gauge transforming the temporal
solutions. Property (iii), essentially described in Section~\ref{subsec:temp-small}, 
allows us to reverse the process, and show
that small energy global well-posedness in the caloric gauge implies
small energy global well-posedness in the temporal gauge.  Finally, as
a consequence of property (ii) the small energy global well-posedness
in the temporal gauge implies large energy local well-posedness in the
temporal gauge. Even better, it shows that the local solutions can be
continued in the temporal gauge for as long as no energy concentration 
occurs in a light cone.

\subsubsection{ Finite energy solutions}

A consequence of \cite{KT} and of the first two papers in the series
\cite{OTYM1}, \cite{OTYM2} is that the small data problem for the
$4+1$ dimensional hyperbolic Yang--Mills equation is well-posed in
several gauges: Coulomb, caloric, and temporal.  In \cite{OTYM2.5} we
exploit the temporal gauge small data result, combined with causality,
to obtain results for the large data problem.  The local in time
result is as follows:

\begin{theorem}[\cite{OTYM2.5}]\label{t:local-temp}
\begin{enumerate}
\item For each finite energy data $(a,e)$ in $\bbR^4$  with concentration scale $r_c$ 
there exists a unique finite energy  solution $A$ to \eqref{ym} in the time interval $[-r_c,r_c]$ in the 
temporal gauge $A_0 = 0$, depending continuously on the initial data. Furthermore, 
any other finite energy solution with the same data must be gauge equivalent to $A$.
\item The same result holds for data in a ball $B_R$ and  the solution in
the corresponding domain of uniqueness $\calD(C_R) \cap ([-r_{c}, r_{c}] \times \bbR^{4})$.
\end{enumerate}
\end{theorem}

We remark that this caloric gauge well-posedness result is in some sense a soft result, which 
is not accompanied by any dispersive type estimates. In expanded form, it asserts that 
regular data generates regular solutions on the $r_c$ time scale, and that the data to solution 
map has a continuous extension to all finite energy data in the uniform energy norm.
However, its proof is anything but straightforward, as it requires the full strength of the 
local well-posedness in the caloric gauge.

Now we consider the continuation question. The next result asserts that temporal solutions can be continued
until energy concentration (i.e. blow-up) occurs. Thus, temporal solutions are also maximal solutions
for the Yang--Mills equation.

\begin{theorem}\label{t:global-temp}
\begin{enumerate}
\item For each finite energy data $(a,e)$ in $\bbR^4$, let $(T_{min},T_{max})$ be the maximal time interval on which the 
temporal gauge solution exists. If $T_{max}$ is finite then we have
\[
\lim_{t \to T_{max}} r_c(t) = 0
\]
and similarly for $T_{min}$. Furthermore, there exits some $X \in \R^4$ so that energy concentrates 
in the backward light cone of $(T_{max},X)$ (respectively the forward light cone of $(T_{min},X)$).

\item The same result holds for data in a ball $B_R$ and  the solution in
the corresponding domain of uniqueness $\calD(B_R)$.
\end{enumerate}
\end{theorem}
The main advantage of this theorem is that it allows us to work with
solutions which do not admit a global caloric representation. The vanishing of $r_c(t)$ is a corollary of Theorem~\ref{t:local-temp}, while existence of a energy concentration point follows by a standard argument; see, e.g., \cite[Lemma~8.1]{OT3}.

The temporal gauge is convenient in order to deal with causality, but not so much in terms of regularity, as it
lacks good $S^{1}$ bounds. For this reason it is convenient to borrow the caloric gauge regularity:

\begin{theorem} \label{t:good-gauge}
  Let $A$ be a finite energy Yang--Mills solution in a cone section
  $C_{[t_1,t_2]}$ with energy concentration scale $r_c$.  Then in
  a suitable gauge $A$ satisfies the bound
\begin{equation}\label{good-gauge}
\| A\|_{L^\infty (\dot H^1 \cap L^4)} + \| \partial_t A\|_{L^\infty L^2}+ 
\| \partial^j A_j \|_{\ell^1 \dot H^\frac12} + \|A_0\|_{\ell^1 \dot H^{\frac32}} 
+ \| \Box A_x\|_{L^2 \dot H^{-\frac12}} \lesssim_{E,\frac{r_c}{t_2}} 1
\end{equation}
in the smaller cone $C_{[t_1,t_2]}^{4r_c}$ where the radius has been decreased by $4 r_c$.
\end{theorem}
The proof of this theorem requires a good gluing technique for local connections
with suitable regularity, which were used to prove Proposition~\ref{prop:goodrep-ball} and Theorem~\ref{thm:goodrep} as well.

\section{To bubble or not to bubble} \label{subsec:no-bubble}

In this section we outline the proof of our two main results in
Theorems~\ref{thm:threshold} and \ref{thm:dichotomy}, following our fourth and the final article \cite{OTYM3}. 
This is based on a blow-up argument  based on Morawetz-type monotonicity formulas,
broadly following the outline of prior works on Wave Maps \cite{ST2} and Maxwell--Klein--Gordon \cite{OT3}.
However, new difficulties arise here both at the conceptual level and at the technical level
due to the more nonlinear gauge features inherent in Yang--Mills and to the nontrivial topological 
structure.

We start with a common part to both proofs, namely a energy-based criterion for soliton bubbling-off, and then we consider the two results separately.

\subsection{A bubble-off criterion}

Our aim here is to describe the proof of the following result, which provides a bubbling-off criterion that applies 
equally for both the Threshold and the Dichotomy Theorems.

\begin{theorem}[Bubbling Theorem]\label{thm:bubble-off}
\begin{enumerate}
\item Let $A$ be a finite energy Yang--Mills wave which blows up in finite time at $(T, X)$.
Assume in addition that for some $\gamma < 1$ we have
\begin{equation} 
\limsup_{t \nearrow T} \calE_{C_{\gamma} \cap S_{t}}(A) > 0, \qquad C_{\gmm} = \{ |x - X| \leq \gmm \abs{t - T }\}.
\end{equation}
Then $A$ bubbles off a soliton at $(T, X)$, as described after Theorem~\ref{thm:dichotomy}.

\item Let $A$ be a finite energy Yang--Mills wave which is global forward in time.
Assume in addition that for some $\gamma < 1$ we have
\begin{equation} 
\limsup_{t \nearrow \infty} \calE_{C_{\gamma} \cap S_{t}}(A) > 0, \qquad C_{\gmm} = \{ |x| \leq \gmm t\}.
\end{equation}
Then $A$ bubbles off a soliton at infinity, as described after Theorem~\ref{thm:dichotomy}.
\end{enumerate}
\end{theorem}

\subsubsection*{Beginning of the proof}
We start with some notations and initial simplifications. 
In the finite time blow up case,
by translation and reflection we can assume that $(T,X) = (0,0)$, and
that the blow-up occurs in the forward light cone.  We introduce the
forward cone $C$, its lateral boundary $\rd C$ and the foliation
$\set{S_{t}}_{t \in [0, \infty)}$ as
\begin{equation*}
	C = \set{(t, x) : 0 \leq \abs{x} \leq t}, \quad
	\rd C = \set{(t, x) : 0 \leq \abs{x} = t}, \quad
	S_{t} = C \cap (\set{t} \times \bbR^{4}).
\end{equation*}
We introduce the energy flux $\calF_{[t_{1}, t_{2}]}(A)$, defined as
\begin{equation*}
	\calF_{[t_{1}, t_{2}]}(A) = \calE_{t_{2}}(A) - \calE_{t_{1}}(A).
\end{equation*}
Assume, for simplicity, that $A$ is regular. Then in both scenarios, by the above energy flux relation, 
we can easily obtain a sequence $A^{(n)}$ of Yang--Mills waves, all obtained by rescaling
the original $A$, and having the following properties:

\begin{enumerate}
\item $A^{(n)}$ is defined on $C_{[\veps_{n}, 1]}$ where $\veps_{n} \to 0$;
\item (Bounded energy in the cone) $\nE_{S_{t}}(A^{(n)}) \leq E$ for every $t \in [\veps_{n}, 1]$;
\item (Decaying flux on $\rd C$) $\calF_{[\veps_{n}, 1]}(A^{(n)}) \leq \veps_{n}^{\frac{1}{2}} E$;
\item (Nontrivial time-like energy at $t = 1$) $\nE_{C_{\gmm} \cap S_{1}}(A^{(n)}) \geq E_0 > 0 $.
\end{enumerate}

\subsubsection*{A Morawetz identity}
Here we describe the key monotonicity formula (or a Morawetz identity), from which we obtain both asymptotic stationarity and compactness for bubble extraction. The idea is to use the \emph{renormalized scaling vector field} $X_{0} = \frac{1}{\sqrt{t^{2} - \abs{x}^{2}}}(t \rd_{t} + x \cdot \rd_{x})$ as a multiplier. Introducing
\begin{equation*}
	{}^{(X_{0})} P_{\alp} (A) = T_{\alp \bt} (A) X_{0}^{\bt} ,
\end{equation*}
where $T_{\alp \bt}(A)$ is the Yang--Mills energy-momentum tensor, we have
\begin{equation} \label{eq:morawetz-0}
	\mathrm{div} \, {}^{(X_{0})} P (A)= \frac{2}{\rho_{0}} \abs{\iota_{X_{0}} F}^{2},
\end{equation}
where $\rho_{0} = \sqrt{t^{2} - \abs{x}^{2}}$. Remarkably, the RHS is nonnegative! 

To derive a monotonicity formula, we would like to integrate \eqref{eq:morawetz-0} on $C_{[t_{1}, t_{2}]}$ and apply the divergence theorem. However, this is not possible since the weight $\rho^{-1}$ blows up on $\rd C$. Instead we introduce a parameter $\veps > 0$ and consider  $X_{\veps} = \frac{1}{\rho_{\veps}} ((t+\veps) \rd_{t} + x \cdot \rd_{x})$, where $\rho_{\veps} = \sqrt{(t+\veps)^{2} - \abs{x}^{2}}$.
Introducing the notation
\begin{equation*}
	{}^{(X_{\veps})} \calP_{S_{t}}(A) = \int_{S_{t}} {}^{(X_{0})} P_{0} (A) \, \ud x,
\end{equation*}
we arrive at 
\begin{equation} \label{eq:int-mono}
\begin{aligned}
	{}^{(X_{\veps})} \calP_{S_{t_{2}}}(A)
	+ \int_{C_{[t_{1}, t_{2}]}} \frac{1}{\rho_{\eps}} \abs{\iota_{X_{\eps}} F}^{2} \, \ud t \ud x 
	= & {}^{(X_{\veps})} \calP_{S_{t_{1}}} + \int_{\rd C_{[t_{1}, t_{2}]}} {}^{(X_{\eps})}P_{\alp} (A) L^{\alp} \, \ud \mathrm{Area}
\end{aligned}
\end{equation}
where $L = \rd_{t} + \frac{x}{\abs{x}} \cdot \rd_{x}$. In the ideal case when the integral on $\rd C$ vanishes, \eqref{eq:int-mono} says that the quantity ${}^{(X_{\veps})} \calP_{S_{t}}$ is monotone in $t$. 

To describe ${}^{(X_{\veps})} \calP_{S_{t}}$ in detail, we need more notation. Let $L = \rd_{t} + \frac{x}{\abs{x}} \cdot \rd_{x}$, $\underline{L} = \rd_{t} - \frac{x}{\abs{x}} \cdot \rd_{x}$, and let $\set{e_{\mathfrak{a}}}_{\mathfrak{2, 3, 4}}$ be orthonormal vectors which are orthogonal to $L, \underline{L}$. In terms of the null decomposition of $F$ defined as
\begin{equation*}
	\alp_{\mathfrak{a}} = F(L, e_{\mathfrak{a}}), \quad
	\underline{\alp}_{\mathfrak{a}} = F(\underline{L}, e_{\mathfrak{a}}), \quad
	\varrho = \frac{1}{2} F(L, \underline{L}), \quad
	\sgm_{\mathfrak{a} \mathfrak{b}} = F(e_{\mathfrak{a}}, e_{\mathfrak{b}}),
\end{equation*}
we have
\begin{equation} \label{eq:mono-eps}
\begin{aligned}
	{}^{(X_{\veps})} \calP_{S_{t}}(A)
	=& \int_{S_{t}} \bigg( \frac{1}{2} \left(\frac{t+r + \veps}{t-r + \veps}\right)^{1/2} \left( \abs{\alp}^{2} + \abs{\varrho}^{2} + \abs{\sgm}^{2} \right) \\
	& \phantom{\int_{S_{t}} \bigg(} 
	+ \frac{1}{2} \left(\frac{t-r + \veps}{t+r + \veps}\right)^{1/2} \left( \abs{\underline{\alp}}^{2} + \abs{\varrho}^{2} + \abs{\sgm}^{2} \right) \bigg) \, \ud x.
\end{aligned}
\end{equation}

Finally, we discuss how \eqref{eq:int-mono} is applied to our setting.
For the solution $A^{(n)}$ constructed above, the RHS of \eqref{eq:int-mono} can be bounded by $\aleq E$ for $\veps = \veps_{n}$. We point out that the last term is bounded by the energy flux $\calF_{[t_{1}, t_{2}]}(A)$. Thus
\begin{equation} \label{eq:X-est}
	\sup_{t \in (\veps_{n}, 1]} {}^{(X_{\veps_{n}})}\calP_{S_{t}}(A^{(n)}) + \iint_{C_{(\veps_{n}, 1]}} \frac{1}{\rho_{\veps_{n}}} \abs{\iota_{X_{\veps_{n}}} F^{(n)}}^{2} \, \ud t \ud x \aleq E.
\end{equation}

Consider a time-like cone $C_{\gmm} = \set{(t, x) : \abs{x} \leq \gmm t}$ for any $0 < \gmm < 1$. Observe that $\rho_{\veps} \simeq t$ and $X_{\veps}$ is uniformly time-like in $C_{\gmm} \cap \set{t \geq 2 \veps}$ (both statements are uniform as $\veps \to 0$ but degenerate as $\gmm \to 1$). Thus boundedness of the spacetime integral term in \eqref{eq:X-est} implies logarithmic integrated decay of a uniformly time-like interior derivative of $F^{(n)}$ in $C_{\gmm}$; this decay is the source of asymptotic stationarity and compactness.

\subsubsection*{Propagating energy in time-like region}
The monotonicity formula \eqref{eq:int-mono} suggests that the weighted energy ${}^{(X_{0})} \calP_{S_{t}}(A^{(n)})$ essentially increases toward the tip.
Using a suitably localized version of the formula, we show that nontrivial energy persists in a time-like cone toward the tip:
\begin{equation} \label{eq:tlike-en}
	\nE_{C_{\gmm} \cap S_{t}} (A^{(n)}) \geq E_{1} \quad \hbox{ for } t \in [\veps_{n}^{\frac{1}{2}}, \veps_{n}^{\frac{1}{4}}],
\end{equation}
where we make $1-\gmm$ and $E_{1}$  smaller if necessary.

\subsubsection*{Final rescaling}
After a pigeonhole argument and suitable rescalings, we obtain a sequence of caloric Yang--Mills waves on $[1, T_{n}] \times \bbR^{4}$ (where $T_{n} \to \infty$), which we still denote by $A^{n}$, with the following properties (final rescaled sequence):
\begin{enumerate}
\item (Bounded energy in the cone) $\nE_{S_{t}}(A^{(n)}) \leq E \qquad (t \in [1, T_{n}]$);
\item (Small energy outside the cone) $\nE_{(\set{t} \times \bbR^{4}) \setminus S_{t}} (A^{(n)}) \ll E \qquad (t \in [1, T_{n}]$);
\item (Nontrivial energy in a time-like region) $\nE_{C_{\gmm} \cap S_{t}} (A^{(n)}) \geq E_{1} \qquad (t \in [1, T_{n}]$);
\item (Asymptotic self-similarity) For every compact subset $\tilde{C}$ of $C_{[1, \infty)}^{1} = \set{(t,x) \in C : \abs{x - \abs{t}} \geq 1}$,
\begin{equation} \label{eq:asymp-ss}
	\iint_{\tilde{C}} \abs{\iota_{X_{0}} F^{(n)}} \, \ud t \ud x \to 0 \quad \hbox{ as } n \to \infty.
\end{equation}
\end{enumerate}

\subsubsection*{Locating concentration scales}
To extract a bubble, we now locate (locally) smallest concentration scales in $A^{(n)}$, which retains the decay \eqref{eq:asymp-ss}. 
A combinatorial argument from \cite{OT3} (based on \cite{ST2}) establishes two possible scenarios (along a subsequence of $A^{(n)}$):
\begin{enumerate}[label=\roman*)]
\item (Time-like concentration) There exists $r > 0$, a sequence of points $(t_{n}, x_{n}) \to (t_{0}, x_{0}) \in \mathrm{Int}(C_{[1, \infty)})$, and a sequence of scales $r_{n} \to 0$ such that
\begin{equation*}
\sup_{x \in B_{r}(x_{n})} \nE_{B_{r_{n}}(x)}(A^{(n)})
\end{equation*}
is uniformly small but nontrivial, yet
\begin{equation*}
	\frac{1}{2r_{n}} \int_{t_{n}-r_{n}}^{t_{n}+r_{n}} \int_{B_{r}(x_{n})} \abs{\iota_{V} F^{(n)}} \, \ud t \ud x \to 0 \quad \hbox{ as } n \to \infty.
\end{equation*}
 where $V = X_{0}(t_{0}, x_{0})$.
\item (Self-similar concentration) For every set\footnote{In fact, any compact subset $\tilde{C}$ in the interior of $C_{[1, \infty)}^{1}$ would work} of the form
\begin{equation*}
\tilde{C} = \set{(t, x) : 0 \leq \abs{x} < t - \frac{1}{2}, \, 2^{j} \leq t < 2^{j+1} \hbox{ for some } j \in \bbZ}
\end{equation*}
there exists $r = r(\tilde{C})$ such that 
\begin{equation*}
\sup_{x \in \tilde{C}} \nE_{B_{r}(x)}(A^{(n)})
\end{equation*}
is uniformly small.
\end{enumerate}

\subsubsection*{Local compactness result}
In both scenarios, we would like to extract a limit modulo scalings,
translations and gauge transformations.  To ensure that the limit is
nontrivial and solves the hyperbolic Yang--Mills equation, we need a
means to ensure compactness.
\begin{theorem} \label{thm:cpt} 
Let $A^{(n)}$ be a sequence of finite energy Yang--Mills connections in $[-2,2] \times \R^4$
which is locally uniformly bounded in the sense of \eqref{good-gauge}.
 Let $Q = [-1, 1]  \times B_{R}(0)$ and $2Q = [-2, 2] \times B_{2R}(0)$. Assume that
  \begin{equation*}
    \lim_{n \to \infty} \nrm{\iota_{X} F}_{L^{2}(2 Q)} = 0,
  \end{equation*}
  where $X$ is a smooth time-like vector field.
  Then on a subsequence, we have
  \begin{equation*}
    A^{(n)} \to A \quad \hbox{ in } H^{1}(Q),
  \end{equation*}
  where $A$ is a solution to the Yang--Mills equation satisfying
  $\iota_{X} F = 0$.
\end{theorem}

The idea of the proof is as follows. The $S^{1}$ bound implies uniform
boundedness of $\nrm{\Box A^{(n)}}_{L^{2}
  \dot{H}^{-\frac{1}{2}}}$. This in turn implies extra regularity away
from the characteristic cone $\set{\abs{\tau} = \abs{\xi}}$ in
frequency space, since $\Box$ is elliptic there. Near the
characteristic cone, we use the following equation for $A^{(n)}$:
\begin{align*}
  X^{\alp} \rd_{\alp} A^{(n)}_{j} - X^{\ell} \rd_{j} A^{(n)}_{\ell} =& -(\iota_{X} F^{(n)})_{j} + \hbox{(smoother error)}, \\
  X^{\ell} \rd_{0} A^{(n)}_{\ell} =& - (\iota_{X} F^{(n)})_{0} +
  \hbox{(smoother error)} .
\end{align*}
Although the system on the LHS is not elliptic, it is microlocally
elliptic (of order $1$) near the characteristic cone $\set{\abs{\tau}
  = \abs{\xi}}$ in frequency space. Inverting this system, and using
the hypothesis $\iota_{X} F^{(n)} \to 0$ in $L^{2}(2 Q)$, we arrive at the
decomposition
\begin{equation*}
  A^{(n)} = A^{(n),small} + A^{(n),smooth}, \quad \nrm{A^{(n),small}}_{H^{1}(Q)} \to 0, \quad \nrm{A^{(n),smooth}}_{H^{1+\alp}(Q)} \aleq 1,
\end{equation*}
for some $\alp > 0$ (in fact, $\alp = \frac{1}{2}$). Applying
Rellich--Kondrachov to $A^{(n), smooth}$, the theorem follows.

\subsubsection*{Extraction of limiting profiles}
In order to apply Theorem~\ref{thm:cpt} in Scenario i), we rescale and 
translate so that $B_{r_{n}}(x_{n}) \to B_{1}(0)$ and apply Theorem~\ref{t:good-gauge} 
to insure the bound \eqref{good-gauge}, uniformly on bounded sets. 
As a result, we extract a nontrivial finite energy
stationary solution (i.e., a soliton).

In Scenario ii), we apply a similar procedure to $B_{r}(0)$, where we
rely on Property~(4) of the final rescaled sequence for the decay
hypothesis in Theorem~\ref{thm:cpt}. In this case, we extract a finite
energy self-similar solution on $C_{[1, \infty)}^{1}$, which is nontrivial
thanks to Property~(3).

\subsubsection*{Exclusion of the self-similar case}
To conclude the bubble extraction argument, it remains to rule out Scenario~ii), i.e., to prove that every finite energy self-similar solution is trivial. 

By self-similarity, the solution restricted to the hyperbolic space $\bbH^{4} = \set{(t, x) : t > 0, \, t^{2} - \abs{x}^{2} = 1}$ is a harmonic Yang--Mills connection. Recall that the harmonic Yang--Mills equation in dimension $4$ is conformally invariant. Thus, by a stereographic projection, we obtain a harmonic Yang--Mills connection on $\bbD^{4}$, which we still denote by $A$. Finite energy condition restricted to the hyperbolic space $\bbH^{4}$ essentially implies that, after a suitable gauge transformation, $A$ is smooth up to the boundary and $A \restriction_{\rd \bbD^{4}}$ vanishes. By an elliptic unique continuation argument (applied to $F$), it follows that the solution is trivial.

\subsection{The Threshold Theorem}

We first  restate our Threshold Theorem  in the caloric gauge.
We will consider the global solvability question for the system \eqref{ym}
 with initial data at time $t = 0$
\begin{equation}\label{data}
(A_j(0), \partial_0 A_j(0)) = (A_{0j}, B_{0j}) \in T^{L^2} \calC \subset 
 \H:= \bfH(\R^4) \times L^2(\R^4).
\end{equation}
Here the caloric gauge imposes a constraint on both $A_{0j}$ and on
$B_{0j}$.  As discussed before, the temporal components of the
connection, namely $A_0$ and $\partial_0 A_0$, are determined in an
elliptic fashion in terms of $A_x$ and $\partial_0 A_x$.

We will also consider higher regularity and (weak) Lipschitz dependence properties of the solutions, using the
spaces 
\[
\H^{\sgm} = \dot{\H}^{\sgm} \cap \calH, \quad
\dot{\H}^{\sgm} = \dot{H}^{\sgm}(\R^{4}) \times \dot{H}^{\sgm-1}(\R^{4}).
\]
Now we can provide a more complete statement for  our main result:

\begin{theorem}\label{t:main}
The Yang--Mills system in the caloric  gauge \eqref{ym}
is globally well-posed in $\H$  for  all caloric initial data  in $\H$
below the ground state energy, in the following sense:

(i) (Regular data) If in addition the data set $(A_{0j}, B_{0j})$ is more
regular, $(A_{0j}, B_{0j}) \in \H^N$, then there exists a unique
global regular caloric solution $(A_j,\partial_0 A_j) \in C(\R, \H^N)$, also
with $(A_0,\partial_0 A_0) \in C(\R, \H^N)$, which has Lipschitz
dependence on the initial data locally in time in the $\H^N$ topology.

(ii) (Rough data) The flow map admits an extension
\[
T^{L^2} \calC  \ni (A_{0j}, B_{0j}) \to (A_\alpha, \partial_t A_\alpha) \in C(\R, \H)
\]
and which is continuous
in the $\H \cap \dot \H^s$ topology for $s <1$ and close to $1$.

(iii) (Weak Lipschitz dependence) The flow map is globally Lipschitz in the 
$\dot \H^s$ topology for $s < 1$, close to $1$.
\end{theorem}

We remark that in effect the proof of the theorem provides a stronger
statement, where the regularity of the solutions is described in terms
of function spaces $S^1$, $S^N$ which incorporate, in particular, Strichartz norms,
$X^{s,b}$ norms and null frame spaces.

Implicit in Theorem ~\ref{t:main} is also a scattering result;
however, this is not so easy to state as it is a modified rather than
linear scattering. In a weaker sense, one can think of scattering as
simply the fact that the $S^1$ norm is finite.

In what follows we outline the proof, using Theorems~\ref{thm:lwp-simple}, \ref{thm:ed-simple} and \ref{thm:bubble-off} 
as our starting point.

\subsubsection{ No bubbling} The first step here is to show that no bubbling
can occur. Here, we closely follow the argument in \cite{LO}.

Indeed, assume by contradiction that a sequence $A^{(n)}$ of rescales and translates
of $A$ converges locally in $H^1$ to a Lorentz transform of a nontrivial soliton $L_v Q$, which implies $L^{2}_{loc}$ convergence of curvature tensors $F^{(n)}$. So after taking a subsequence, for almost every $t$
\begin{equation*}
	\calE_{\set{t} \times B_{R}}(A^{(n)}) = \frac{1}{2} \int_{B_{R}} \brk{F^{(n)}, F^{(n)}}(t) \to \calE_{\set{t} \times B_{R}}(L_{v} Q) \qquad \hbox{ for any } R > 0,
\end{equation*}
which in turn implies
\begin{equation*}
	\nE(Q) \leq \nE(A) < 2 \Egs.
\end{equation*}
By Theorem~\ref{thm:gs}, the only possibility for $Q$ is that $\abs{\ch(Q)} = \spE(Q)$. Moreover, since Lorentz transforms preserve the topological class, $\ch(L_{v}(Q)) = \ch(Q)$. 

By topological triviality of $A^{(n)}(t)$, we have $\ch(A^{(n)}(t)) = 0$, and thus
\begin{equation*}
	\int_{\bbR^{4} \setminus B_{R}(0)} - \brk{F^{(n)} \wedge F^{(n)}}(t) = - \int_{B_{R}(0)} - \brk{F^{(n)} \wedge F^{(n)}}(t).
\end{equation*}
By $L^{2}_{loc}$ convergence of $F^{(n)}$, the absolute value of the first term on the RHS can be made arbitrarily close to $\abs{\ch(A)} = \nE(Q)$ by taking $R$ very large. Using the Bogomoln'yi lower bound $\abs{\brk{F \wedge F}} \leq \frac{1}{2} \brk{F_{ij}, F^{ij}}$ in $\bbR^{4} \setminus B_{R}$, it follows that
\begin{align*}
	\nE(A) 
	\geq & \limsup_{n \to \infty} \Big( \frac{1}{2} \int_{B_{R}} \brk{F^{(n)}, F^{(n)}}(t) + \abs{\int_{\bbR^{4} \setminus B_{R}} \brk{F^{(n)} \wedge F^{(n)}}(t)} \Big) \\
	\geq & \nE_{\set{t} \times B_{R}}(L_{v} Q) + \abs{\int_{B_{R}} -\brk{F[L_{v} Q] \wedge F[L_{v} Q]}} \\
	\geq & \nE(L_{v} Q) + \nE(Q) - o_{R \to \infty}(1) .
\end{align*}
Since $\nE(L_{v} Q) \geq \nE(Q) \geq \Egs$, we reach a contradiction.

\subsubsection{No blow-up}
Suppose finite time blow-up occurs for a subthreshold caloric Yang--Mills
wave. By translation invariance we can assume that the blow-up happens at $(0,0)$,
backwards in time. By the small data result, we must have energy concentration in the 
forward light cone $C$ at $t = 0$
\begin{equation}\label{en-conc}
\lim_{t \searrow 0} \nE_{S_{t}}(A) > 0 .
\end{equation}
On the other hand, as bubbling cannot occur, by Theorem~\ref{thm:bubble-off} we must have
\begin{equation}\label{small-in}
\lim_{t \searrow 0} \nE_{C_{\gmm} \cap S_{t}}(A) = 0 \qquad \forall \gamma < 1 .
\end{equation}

To reach a contradiction, it would suffice to show that the energy dispersion  decays 
near the tip of the cone,
\[
\lim_{t \searrow 0} \| F\|_{ED[0,t]} = 0 .
\]
Then Theorem~\ref{thm:ed-simple} would yield a bound for $\| A\|_{S^1[0,t]}$,
which shows that the solution $A$ extends below $t = 0$ and in
particular the energy concentration \eqref{en-conc} cannot occur.

One problem with this strategy is that we have no a-priori knowledge about what happens
outside the cone. To rectify this we excise the outer part of the solution, so that we 
are left with a connection $\tA$ in a small time interval $[0,t_0]$, so that 
\begin{enumerate}
\item The two connections agree inside, $\tA = A$ in $C_{[0,t_0]}$.
\item $\tA$ has small energy outside, 
\begin{equation}\label{small-out}
\nE_{\R^4 \setminus C_t}(\tA) \leq \epsilon \ll 1, \qquad t \in [0,t_0]
\end{equation}
\end{enumerate}
Here $\epsilon$ can be chosen arbitrarily small, and $t_0$ depends on $\epsilon$.
This is achieved using Proposition~\ref{prop:id-excise} at a well
chosen time $t_0$, using the flux decay near the tip of the cone.  By
finite speed of propagation, note that the new and old solutions agree
in $C$.  In particular, the new solution also concentrates energy at
$(0, 0)$, and thus cannot be extended past $0$.

Taking into account \eqref{small-out} and \eqref{small-in} (the latter transfers from $A$ to $\tA$)
for $\tA$, we see that the energy of $\tA$ must concentrate near the cone. Using 
the Morawetz estimate \eqref{eq:X-est}, we obtain as well a second energy bound inside the cone,
namely 
\begin{equation}\label{X-gamma}
\limsup_{t \to 0} {}^{(X^\gamma)} \calP_{S_{t}}[\tA] \lesssim_\nE 1, \qquad \gamma < 1.
\end{equation}
This shows that in addition, only certain curvature components may be large near the cone.

Finally, we are in a position to show that $\tA$ is energy dispersed near the tip, and 
thus reach the desired contradiction by Theorem~\ref{thm:ed-simple}. This is done using the 
following result:

\begin{proposition}\label{p:small-ed}
 Let $(A_{x}, \rd_{t} A_{x})(t)$ be a caloric Yang--Mills data with energy $\nE < 2\Egs$. Then for each $\epsilon > 0$
there exists $\gamma < 1$ and $\delta > 0$ so that the bounds 
\[
\nE_{C_{\gmm} \cap S_{t}(A)}(A) + \nE_{\R^4 \setminus S_{t}}(\tA) \leq \delta, \qquad  {}^{(X^\gamma)} \calP_{S_{t}}[\tA] \lesssim_\nE 1
\]
imply 
\[
 \| F\|_{ED[t]}  \leq \epsilon.
\]
\end{proposition}

Indeed, by the huge weight near $\rd C$ in
${}^{(X_{\veps_{n}})} \calP_{S_{1}}(A)$ and smallness of energy
elsewhere, all components of $F$ except for $\underline{\alp}$ are
small in $L^{2}$. To control $\underline{\alp}$, it suffices to
consider $F_{r \mathfrak{a}} = \alp_{\mathfrak{a}} -
\underline{\alp}_{\mathfrak{a}}$ in the frame $(e_{t} = \rd_{t}, e_{r}
= \rd_{r}, e_{2}, e_{3}, e_{4})$. By the Yang--Mills equation and the
Bianchi identity, they obey the following covariant div-curl system on
spheres\footnote{We remark that in our actual proof, we work with an
  analogous div-curl system on hyperplanes for technical simplicity.}:
\begin{align*}
  \covD_{\mathfrak{a}} F_{r \mathfrak{b}} - \covD_{\mathfrak{b}} F_{r \mathfrak{a}} =& \covD_{r} \sgm_{\mathfrak{a} \mathfrak{b}}, \\
  \covD^{\mathfrak{a}} F_{r \mathfrak{a}} =& \covD^{\mathfrak{a}}
  \alp_{\mathfrak{a}} + \covD_{r} \varrho.
\end{align*}
The crucial observation is that the RHS only involve components with
small energy. In the abelian case (where $\covD = \nb$), this div-curl
system can be easily inverted, and it follows that
$\nrm{\abs{\nb_{x}}^{-1} \anb F_{r \mathfrak{a}}}_{L^{2}} \ll E$,
where $\anb = (\nb_{e_{2}}, \nb_{e_{3}}, \nb_{e_{4}})$ stands for the
angular derivatives. By Bernstein, this is sufficient to rule out the
null concentration scenario. A more involved  argument is needed 
in the non-abelian case.

\subsubsection{Scattering}
The argument here is similar but simpler. Simply by translating the coordinate system 
we can insure that the condition \eqref{small-out} holds for $t \in [t_0,\infty)$.
Then the rest of the argument carries through unchanged.

\subsection{The Dichotomy Theorem}

Here we would like to apply the same argument as before. This time we
are assuming, rather than proving that bubbling does not happen. We
can still truncate the solution $A$ outside to insure that the bound
\eqref{small-out} holds in the blow-up case, or translate the
coordinates to achieve the same outcome in the non-scattering case.
The new difficulty is that we are no longer guaranteed that we can work 
in the caloric gauge, as the energy may be above the threshold.

However, it turns out that this is only a technical obstruction, as we 
can now prove a much stronger form of 

\begin{proposition}\label{p:small-ed+}
 Let $(A_{x}, F_{0x})(t)$ be a finite energy Yang--Mills data with energy $\nE$. Then for each $\epsilon > 0$
there exists $\gamma < 1$ and $\delta > 0$ so that the bounds 
\[
\nE_{C_{\gmm} \cap S_{t}(A)}(A) + \nE_{\R^4 \setminus S_{t}}(\tA) \leq \delta, \qquad  {}^{(X^\gamma)} \calP_{S_{t}}[\tA] \lesssim_\nE 1
\]
imply that $(A_{x}, F_{0x})(t)$ admits a caloric gauge representation so that in addition we have
\[
 \| F\|_{ED[t]}  \leq \epsilon
\]
\end{proposition}

The difficulty here is to obtain the caloric gauge representation,
without assuming any a-priori bound on $\|A[t]\|_{\dot H^1 \times L^2}$. 
This is done via multiple continuity arguments, in several steps:

\begin{enumerate}[label=(\roman*)]
  \item Working in an annulus, use a continuity argument show that one
    can obtain a local gauge which where $A$ is controlled in $\dot
    H^1$, with small $L^4$ norm.

\item Extend previous step to all of $\R^4$, by gluing small $\dot H^1
  \cap L^4$ connections obtained via Uhlenbeck's lemma inside the
  annulus and outside.

\item Use a second continuity argument to show that a corresponding
caloric connection exists.  Here the previous step is used to construct a path to $0$.
\end{enumerate}

\bibliographystyle{ym}
\bibliography{ym}

\end{document}